\documentclass{article}
\usepackage{graphicx} 
\usepackage{amssymb}
\usepackage{color}

\newtheorem{propp}{Proposition}[section]
\newtheorem{theo}[propp]{Theorem}  
\newtheorem{coro}[propp]{Corollary} 
\newtheorem{lemm}[propp]{Lemma} 
\def\proof{\textit{Proof. }}
\newcounter{example}
\newcommand{\example}{\addtocounter{example}{1}\par\noindent\textbf{Example \arabic{example}.} }  

\newcommand{\pic}[4]{\vspace{1ex}\setlength{\unitlength}{1cm}
\begin{picture}(0,#3)(5,.5)
\put(#2){\includegraphics[#4]{#1.eps}}
\end{picture}\vspace{1ex}}

\def\C{{\mathbb C}}
\def\D{{\mathbb D}}
\def\T{{\mathcal T}}
\def\R{{\mathbb R}}
\def\re{\,\!\mbox{Re}\,}
\def\im{\,\!\mbox{Im}\,}
\def\conj{\overline}

\def\qed{\ \hfill\raise3.5pt\hbox{\framebox[1ex]{\ }}}
\newcommand{\rev}[1]{\raisebox{.08ex}{$\stackrel{\leftrightarrow}{#1}$}}

\begin{document}

\begin{center}
  NUMERICAL SOLUTION OF THE BELTRAMI EQUATION\\ 
  VIA A PURELY LINEAR SYSTEM
\end{center} 

\bigskip\noindent
\textrm{R. Michael Porter}%
\footnote{Research partially supported by CONACyT grant 166183}\\
\textrm{Department of Mathematics,\\
Centro de Investigaci\'on y de Estudios Avanzados del I.P.N.,\\
Apdo. Postal 1-798, Arteaga 5,
76000 Quer\'etaro, Qro., Mexico
}\\
\texttt{mike@math.cinvestav.edu.mx} 

\bigskip\noindent
\textrm{Hirokazu Shimauchi}%
\footnote{Research supported by International Advanced Research and Education Organization in Tohoku University }\\
\textrm{Division of Mathematics, \\
Graduate School of Information Sciences, Tohoku University\\
6-3-09 Aramaki-Aza-Aoba, Aoba-ku, Sendai 980-8579, Japan}\\
 \texttt{shimauchi@ims.is.tohoku.ac.jp, hirokazu.shimauchi@gmail.com} 


  
\begin{abstract}
  An effective algorithm is presented for solving the Beltrami
  equation
  $\partial\!f/\overline{\partial}z=\mu\,\partial\!f/\partial\!z$ in a
  planar disk.   The disk is triangulated in a simple way and
  $f$ is approximated by piecewise linear mappings; the images of the
  vertices of the triangles are defined by an overdetermined system of
  linear equations. (Certain apparently nonlinear conditions on the boundary
  are eliminated by means of a symmetry construction.) The linear
  system is sparse and its solution is obtained by standard
  least-squares, so the algorithm involves no evaluation of singular
  integrals nor any iterative procedure for obtaining a single
  approximation of $f$. Numerical examples are provided, including a
  deformation in a Teichm\"uller space of a Fuchsian group.
\end{abstract} 
  
\section{Introduction}
The Beltrami equation
\begin{equation}
  \label{eq:beltrami}
    \frac{ \partial f(z)/ \partial\conj{z} }
         { \partial f(z)/ \partial z } = \mu(z)
\end{equation}
determines a unique normalized quasiconformal self-mapping $f$
of the unit disk $\D=\{z\colon\; |z|<1\}$ in the complex plane.  Here $\mu$
is a given measurable function in $\D$ with $\|\mu\|_\infty<1$, and
is called the Beltrami derivative (or complex dilatation) of $f$.
One says that $f$ is $\mu$-conformal.

The Beltrami equation has been the object of deep investigation in
large measure due to its importance in the theory of deformations of
Kleinian groups and their applications to Teichm\"uller spaces
\cite{Harv,Lehto}.  Other applications of the Beltrami equation as
mentioned in the introduction to \cite{Dar-fastB}  are
quite well known and we will not go into them here.  Some more recent
applications, such as mapping of the cerebral cortex, use the Beltrami
equation in the spirit of its original application, dating back to
Gauss, for finding a conformal mapping from a surface in 3-space onto a
planar region; this is done effectively in \cite{AHTK} although
without explicit use of the Beltrami derivative $\mu$. The Beltrami
derivative has also been proposed as a way of compressing data for
surface maps \cite{LLKWG}.

With this increasing use of computer applications it has become of
great interest to know how solve the Beltrami equation numerically.
One method for doing this is suggested naturally by the classical
existence proof given by Mori-Boyarskii-Ahlfors-Bers
\cite{Harv,Lehto}.  For this method one must evaluate singular
integrals of the form
\[ T_mg(z) = -\frac{1}{\pi}\int\!\!\int_{\D} \frac{g(\zeta)}{(\zeta-z)^m} \,d\xi d\eta
,\quad m=1,2 \]
(defined as Cauchy principal values when $m=2$), and then calculate sums
of Newmann series of the form $\sum T_2(\mu T_2(\cdots(\mu T_2(\mu)\cdots))$.
A related approach involving the singular integrals was 
developed by P.\ Daripa and D.\ Mashat \cite{Dar-Mash,Dar-fastB},
and refined by D.\ Gaidashev and D.\ Khmelev \cite{GK}.  Instead of summing
the Newmann series, their method involves iteration towards a
solution of a related Dirichlet problem. Their work incudes
refinement of the technique of evaluation of the singular integrals
via FFT, which is of interest in itself.  

Zh.-X.\ He \cite{He} proposed an alternative method of solving the
Beltrami equation, based on circle packings. G.\ B.\ Williams \cite{W}
presented another circle-packing method, based on the idea of
conformal welding. Little information is available on the numerical
performance of these methods.  Lui et.\ al \cite{LLSX} describe yet
another method, which reduces the question of solving the Beltrami
equation to that of a linear system on the underlying mesh. Their
method, focused on obtaining Teichm\"uller mappings of prescribed
domains (rather than self-mappings of a disk as we consider here), is
applied to problems of face recognition and brain mapping.  Many other
approaches have been given to solve for quasiconformal mappings, often
from surfaces in $\R^3$ to the plane.  Descriptions and further
references may be found in \cite{GZWLY,LPRM}.  An attempt to solve the
numerical solution of the Beltrami equation by applying conformal
mappings as an intermediate step, was made in \cite{PorArx}. However,
this method was later found not to converge to the proper solution and
appears not to be salvageable.

In this paper we give a much simpler algorithm with full proof of
convergence in the case of a smooth Beltrami derivative. (We believe
that this smoothness condition is overly restrictive, as numerical
experiments suggest.) In Section \ref{sec:prelim} we gather the basic
facts we will need about quasiconformal mappings.  In Sections
\ref{sec:discrete} and \ref{sec:DBE} we explain how to set up the
linear system describing the Beltrami equation.  In Section
\ref{sec:main} the algorithm is specified and the main theorem on the
convergence of the algorithm is stated and proved.  Several numerical
examples are provided in Section \ref{secnumres}, including a
deformation of a Fuchsian group.  In the closing comments we discuss
the corresponding computational cost.

\section{\label{sec:prelim}Preliminaries}

\subsection{Affine linear quasiconformal mappings}\label{sec:affinelinear}

In this section   $\mu$, $a$, $b$ are complex constants subject to
$a\not=0$, $|\mu|<1$, and we consider the mappings
\begin{eqnarray} 
    L_\mu(z) &=& \frac{z + \mu \conj{z}}{1+\mu},  \label{eq:defL}\\
    H_{a,b}(z) &=& az+b  \label{eq:defH},
\end{eqnarray}
for $z\in\C$.  Thus $L_\mu$ is $\mu$-conformal and real-linear, while
$H_{a,b}$ is conformal and affine complex-linear. Note that $L_\mu$ is
determined by its value at any single point other than its fixed points
$z=0$ and $z=1$, while $H_{a,b}$ is determined by the images of any
two points. All $\mu$-conformal affine linear mappings are of the form
$H_{a,b}\circ L_\mu$, and this decomposition is unique.  We will use
the following form of expressing affine linear mappings.

\begin{propp} \label{prop:Bdef} 
  Given $z_1$, $z_2$ distinct and $w_1$, $w_2$ distinct, together with
  $|\mu|<1$, there is a unique $\mu$-conformal affine linear mapping
  $B=B_{\mu;\,z_1,z_2;\,w_1,w_2}$ which sends $z_1$ to $w_1$ and $z_2$
  to $w_2$. This mapping is given explicitly by
\begin{eqnarray} 
 B(z) &=&  \nonumber
     w_1 + \frac{w_2-w_1}{L_\mu(z_2-z_1)}L_\mu(z-z_1)   \\ \nonumber
 &=&       \frac{L_\mu(z_2-z)}{L_\mu(z_2-z_1)}w_1 + \frac{L_\mu(z_1-z)}{L_\mu(z_1-z_2)}w_2 .    \label{bmap}
\end{eqnarray}
\end{propp}
 
The coefficients of $w_1,w_2$ in this last expression are
never $0$, $1$, or $\infty$ when $z_1,z_2,z_3$ are distinct.  As we
will be interested in mappings of triangles, the following facts will
be useful.

\begin{coro} \label{coro:eqntri} 
  If a $\mu$-conformal affine linear map takes $z_1$, $z_2$, $z_3$ to
  $w_1$, $w_2$, $w_3$ respectively, then
\[ L_\mu(z_2-z_3)\,w_1 +  L_\mu(z_3-z_1)\,w_2 +  L_\mu(z_1-z_2)\,w_3 =0.
\] 
\end{coro}

\begin{coro} \label{coro:implicitmu} 
  Given $z_1$, $z_2$, $z_3$ noncollinear and $w_1$, $w_2$, $w_3$
  noncollinear, there is a unique affine linear mapping which sends
  $z_i$ to $w_i$ ($i=1,2,3$).  Its Beltrami derivative is equal to
\begin{equation} \label{eq:implicitmu}
 \mu = -\frac{(z_2-z_1)(w_3-w_1) - (z_3-z_1)(w_2-w_1)}
      {(\conj{z_2}-\conj{z_1})(w_3-w_1) - 
       (\conj{z_3}-\conj{z_1})(w_2-w_1)}.
\end{equation}
\end{coro}

\proof The affine linear mapping is well-defined because each of the
pairs $(z_2-z_1,z_3-z_1)$, $(w_2-w_1,w_3-w_1)$ is linearly independent
over the real numbers.  To calculate the Beltrami derivative, substitute
(\ref{eq:defL}) into the formula of Corollary \ref{coro:eqntri} and solve
for $\mu$. \qed

In applying Corollary \ref{coro:implicitmu}, one normally would also
require the $z$- and $w$-triangles to be like-oriented, to ensure that
$|\mu|<1$.

\subsection{Closeness to similarity}\label{subsec:similarity}

The affine-conformal mapping (\ref{eq:defH}) sends triangles to similar
triangles. A Beltrami derivative (\ref{eq:implicitmu}) can be regarded
as a measure of how much the (positively oriented) triangles with
vertices $z_1,z_2,z_3$ and $w_1,w_2,w_3$ fail to be similar. Indeed,
these triangles may be carried by conformal mappings of the form
$H_{a,b}$ to normalized triangles $(0,1,c_0)$, $(0,1,c)$ respectively,
which in turn correspond under the mapping $L_\mu$ where
\begin{equation}  \label{eq:muparameter}
  \mu = -\frac{c-c_0}{c-\conj{c_0}}.
\end{equation}
With $c_0$ fixed, the relation $c\leftrightarrow\mu$ from $\im
c>0$ to $|\mu|<1$ thus associates an element of the unit disk to each
similarity class of oriented triangles, and we can use $|\mu|$ as a
measure of the discrepancy from being similar to the original
triangle. In this regard we must note that the Beltrami derivative of
the affine mapping from $(z_1,z_2,z_3)$ to $(w_1,w_2,w_3)$ is
$e^{2i\arg(z_2-z_1)}\mu$ (Proposition \ref{prop:bcomposition} below), so
the absolute value is not altered by the normalization to the triangle
$(0,1,c_0)$ with horizontal base.

\subsection{Quasiconformal mappings}
   
The following well known general properties of quasiconformal
mappings \cite{AhlLQM,Lehto,LV} are fundamental to this work.

\begin{propp} \label{prop:existence}
If $\mu$ is measurable in $\D=\{z\in\C\colon|z|<1\}$ and satisfies
$\|\mu\|_\infty<1$, then there is a unique $\mu$-conformal mapping
$f\colon\D\to\D$ satisfying the normalization
\[ f(0)=0,\quad f(1)=1.
\]
\end{propp}
 
\begin{propp} \label{prop:convergence}
 Let $\mu$, $\mu_n$ be measurable functions in $\D$ with $\|\mu_n\|_\infty\le c<1$
and suppose that $\mu_n\to\mu$ pointwise as $n\to\infty$. Let $f$, $f_n$
be the normalized solutions of the corresponding Beltrami equations
given by Proposition \ref{prop:existence}.  Then $f_n$ converge to $f$ uniformly
on compact subsets of $\D$.
 \end{propp}

\begin{propp} \label{prop:bcomposition}
Let $f_1$,$f_2$ be $\mu_1$,$\mu_2$-conformal mappings
respectively.  \\(i) Suppose that $f_1$, $f_2$ are defined in the same
planar domain and
\[   f_2=h\circ f_1. \]
Then $\mu_1=\mu_2$ a.e.\ if and only if $h$ is
a conformal mapping from the image of $f_1$ to the image of $f_2$.
\\(ii) Suppose that   $h$ is defined in 
the domain of $f_2$ and
\[   f_2 =  f_1 \circ h. \]
If $h$ is conformal, then $\mu_2 = (\mu_1\circ h)(\conj{h'}/h')$. 
\\(iii) If 
\[ f_2(z)=\conj{f_1(\conj{z})},\] then $\mu_2(z)=\conj{\mu_1(\conj{z})}$.
\end{propp}

\section{\label{sec:discrete} Context for discrete Beltrami equation}

In this section we describe the geometric and algebraic
elements necessary for our discrete version of the Beltrami equation.

 To discretize the problem, we will consider finite triangulations
 $\T_z$ of the closed unit disk $\overline{\D}$ in the $z$-plane. We always
 assume that the union $|\T_z|$  of the (closed) triangles of $\T_z$ is
 bounded by a Jordan polygon inscribed in $\overline{\D}$.  
  We are given a proposed Beltrami derivative $\mu$ in $\D$, i.e., a
  measurable function such that $\|\mu\|_\infty=\mbox{ess\ sup}_{z \in
  \D}\;|\mu(z)|<1$, and we want to construct an isomorphic (simplicially
  equivalent) triangulation $\T_w$ of the unit disk in the $w$-plane
  such that the induced piecewise linear mapping (PL-mapping) is
  approximately $\mu$-conformal on each triangle $\tau \in \T_z$.
 
\subsection{Logarithmic coordinates}

Given $\mu$, in many cases it is not difficult to find a discrete
$\mu$-conformal mapping to some triangulation $\T_w$ of a domain whose
shape is not predetermined, but the challenge is to make the outer
polygon bounding $\T_w$ precisely a circle.  This question
can be approached many ways, and in general the condition that ``the
points of the outer polygon of $\T_z$ must be mapped to a circle
centered at the origin'' translates into some nonlinear conditions.
In order to evade this ``outer boundary condition'' we can extend
$\mu(z)$ by reflection to the exterior of $\D$ (using the Chain Rule
for Beltrami derivatives), likewise extending $\T_z$ by including the
vertices $1/\conj{z_{jk}}$, and then solve the extended problem on the
Riemann sphere.  However, the nonlinearity of the inversion
$1/\conj{z}$ leads to inaccuracies since it does not respect affine
mappings between triangles. In other words, the inversion in
$\{|z|=1\}$ does not produce a truly symmetric discrete problem. We
prefer to avoid this difficulty as follows, by introducing logarithmic
coordinates in the form of the variables $Z=\log z$, $W=\log w$ in the
left half-plane (cf.\ \cite[section 6.4]{LV}). In other words, in
order to solve for the discrete function $w=f(z)$ from $\D$ to $\D$ we
will solve first for the correspondence
\begin{equation}  \label{eq:WFZ}
   W=F(Z)
\end{equation}
 where $\exp F(Z) = f( \exp z)$, and then simply apply the operation
\begin{equation} \label{eq:expZexpW}
  z=\exp(Z),\quad w=\exp(W)
\end{equation}
to obtain the desired mapping $z\mapsto w$.  Of course, we will need to justify
that the distortion of triangles produced by the exponential mappping does
not affect the accuracy significantly.  

\subsection{Triangulation}
\label{subsec:tri}
\begin{figure}[b!]
  \centering
   \pic{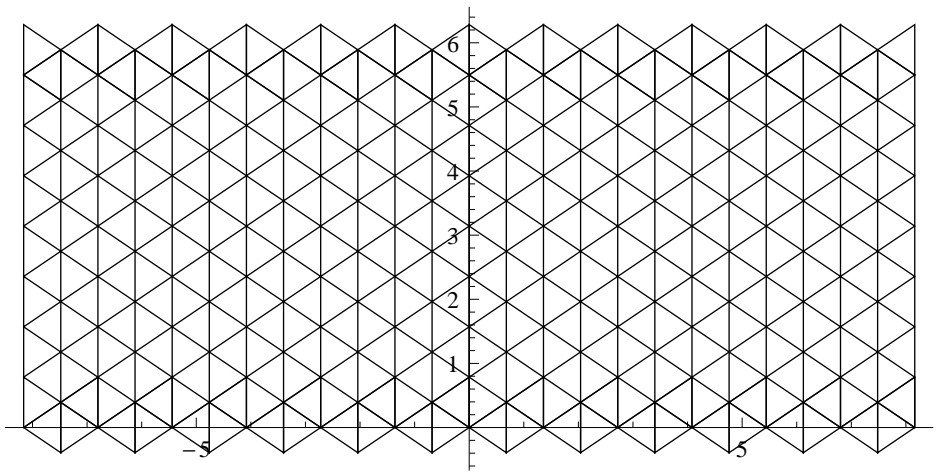}{-1,0}{5.5}{}
  \caption{Basic mesh in $W$-plane, with its reflection in
the imaginary axis.}
  \label{fig:Wmesh}
\end{figure}

It is convenient to use negative indices for points of the basic
domain. We fix the mesh order, that is, a pair of positive integers
$M,N$ (presumably large).  Let
\[ R_{-M}<R_{-M+1}<\cdots<R_{-1}<R_0=0
\]
and define the $(M+1)N$ vertices
\begin{equation}\label{Zjk}
 Z_{jk}= R_j +    \frac{2\pi (k + j_2/2)}{N}i
\end{equation}
for $-M\le j\le 0$ and $0\le k\le N-1$. Here
$j_2=j\bmod 2$ is 0 when $j$ is even and 1 when $j$ is odd.
 
If we extend the formula for arbitrary values of $k$, we obtain
a periodic mesh with $Z_{j,k+N} = Z_{jk}+2\pi i$.
The basic mesh contains $MN$ \emph{rightward pointing triangles}
defined as follows,
\begin{equation}  \label{eq:triright}
  \tau_{jk}^+ = \left\{\begin{array}{ll}
   ( Z_{j-1, k-1}, \ Z_{j-1, k}, \ Z_{j,k}) ,& \mbox{$j$ even,}\\
   ( Z_{j-1, k}, \ Z_{j-1, k+1}, \ Z_{j,k}),& \mbox{$j$ odd,}
   \end{array}\right.
\end{equation}for $-M+1\le j\le0$. There are also
$MN$ \emph{leftward pointing triangles}
\begin{equation}  \label{eq:trileft}
  \tau_{jk}^- =  \left\{\begin{array}{ll}
   (   Z_{j+1, k-1},\  Z_{j+1, k},\ Z_{j, k} ),& \mbox{$j$ even,}\\
   (  Z_{j+1, k},\  Z_{j+1, k+1},\ Z_{j,k}  ),& \mbox{$j$ odd,}
   \end{array}\right.
\end{equation} 
for $-M\le j\le-1$.  Note however, that for even values of $j$ the
triangle $\tau_{j0}^+$ contains $Z_{j-1,-1}$ and the triangle $\tau_{j,0}^-$
contains $Z_{j+1,-1}$; while for odd $j$ the triangle $\tau_{j,N-1}^+$
contains $Z_{j-1,N}$ and the triangle $\tau_{j,N-1}^-$ contains
$Z_{j+1,N}$.  The second index $k$ of each of these points lie outside
of the basic range $0\le k\le N-1$ which most interests us.

We will assume for the rest of this paper that 
\begin{equation}  \label{eq:Rj}
  R_j=(\sqrt{3}\,\pi/N)j,
\end{equation}  
so the triangles $\tau_{jk}^\pm$ will be equilateral. We extend the
structure to the right half-plane using the reflection in the
imaginary axis
\begin{equation}\label{eq:defrho}
 \varrho(Z) = -\conj{Z},
\end{equation}
and setting
\[    Z_{jk} = \varrho(Z_{-j,k})
\]
for indices $j>0$.  Thus in the extended domain we have $-M\le j \le
M$, with symmetry around the index $j=0$.  The extended triangulation
contains $(2M+1)N$ vertices.

Since we are originally given $\mu$ in $\D$, 
we will need the pullback of $\mu$ as a differential of type $(-1,1)$
to the left half-plane  (Proposition \ref{prop:bcomposition}), that is,
\begin{equation} \label{eq:pullback}
 \nu(Z) = \mu(e^Z)\frac{e^{\conj{Z}} }{e^Z} =  \mu(e^Z)e^{-2i\,\im Z}
 ,\quad \re Z<0.
\end{equation}
For $\re Z>0$ we want $F(Z)=\rho(F(\rho(Z))$ according to
(\ref{eq:WFZ}). By Proposition \ref{prop:bcomposition} the Beltrami
derivative of $F$ is $\nu(Z) = \conj{\nu(\varrho(Z))}$.  We will write
$\nu^\pm_{jk}$ for the average value of $\nu(Z)$ on the triangle
$\tau^\pm_{jk}$. While technically this means the integral of $\nu$
divided by the area of the triangle, for numerical work it is
convenient to take the average of $\nu(Z)$ over the three vertices as
an approximation to this integral, at least when $\nu$ is
continuous. Let us note that
\begin{equation}\label{eq:nusymm}
  \nu_{jk} = \conj{\nu_{-j,k}},\quad j>0.  
\end{equation}
    
\subsection{Boundary conditions}
\label{subsec:bdrycond}

We will need to describe the behavior of the discrete $\mu$-conformal
mapping at the points corresponding to $j=\pm M$.  The vertical line
$\{Z\colon\;\re Z = R_{-M}\}$ corresponds to the circle $|z|=r_{-M}$
in the $z$-plane, where we write $r_j=\exp R_j$.  Let $f\colon\D\to\D$
denote the exact (smooth) solution of (\ref{eq:beltrami}), normalized
by $f(0)=0$, $f(1)=1$.  Then near $z=0$, we know (recall
(\ref{eq:defL}), (\ref{eq:defH})) that $f$ is approximately equal to
an affine map $H_{a,0}\circ L_{\mu(0)}$ for some
$a\in\C\setminus\{0\}$, and the image of the small circle is an
ellipse in the $w$-plane whose eccentricity is determined by $\mu(0)$
and whose size is determined by $a$.  Similarly, the right boundary $\{\re Z =
R_{M}\}$ corresponds to the inversion of this small ellipse in the
circumference $\partial\D$ of $\D$. Thus we are considering the
problem as defined in the annulus $1/r_M\le|z|\le r_M$, and the value
$|a|$ is related to the conformal module \cite{LV} of the image of
this annulus under $f$.  We will sidestep the question of explicitly considering
$a$ in \ref{subsec:bdryeq} below.

\section{Discrete Beltrami equation}\label{sec:DBE}
 
Here we define the system of linear equations representing the solution
of the Beltrami equation. First we consider the unknowns 
\[ \{W_{jk}\colon\; -M\le j\le 0,\ 0\le k\le N-1\}
\]
in the left half-plane.

\subsection{Triangle equations}

Following Corollary \ref{coro:eqntri}, to each rightward pointing
triangle $\tau_{jk}^+$ of (\ref{eq:triright}) we associate one linear
equation
 \begin{equation} \label{eq:triangles+}
  a_{jk}^{+}W_{jk} + b_{jk}^{+}W_{j-1,k} + c_{jk}^{+}W_{j-1,k+1} = 0
\end{equation}
where
\begin{eqnarray}
   a_{jk}^{+} &=& \left\{\begin{array}{ll} 
     L_{\nu_{jk}}(Z_{j-1,k-1}-Z_{j-1,k}),&j\mbox{ even},\\ 
     L_{\nu_{jk}}(Z_{j-1,k}-Z_{j-1,k+1}),&j\mbox{ odd},
    \end{array}\right. \nonumber\\ 
   b_{jk}^{+} &=& \left\{\begin{array}{ll} 
     L_{\nu_{jk}}(Z_{j-1,k}-Z_{j,k}),&j\mbox{ even},\\ 
     L_{\nu_{jk}}(Z_{j-1,k+1}-Z_{jk}),&j\mbox{ odd},
    \end{array}\right.  \label{eq:abc+}\\
   c_{jk}^{+} &=&\left\{\begin{array}{ll} 
     L_{\nu_{jk}}(Z_{j,k}-Z_{j-1,k-1}),&j\mbox{ even},\\ 
     L_{\nu_{jk}}(Z_{jk}-Z_{j-1,k}),&j\mbox{ odd}, 
    \end{array}\right.   \nonumber 
\end{eqnarray}
and similarly an equation for each leftward pointing triangle $\tau_{jk}^-$,
\begin{equation} \label{eq:triangles-}
  a_{jk}^{-}W_{jk} + b_{jk}^{-}W_{j+1,k-1} + c_{jk}^{-}W_{j+1,k} = 0
\end{equation}
where
\begin{eqnarray}
   a_{jk}^{-} &=& \left\{\begin{array}{ll}     
     L_{\nu_{jk}}(Z_{j+1,k-1}-Z_{j+1,k}),&j\mbox{ even},\\ 
     L_{\nu_{jk}}(Z_{j+1,k}-Z_{j+1,k+1} ),&j\mbox{ odd},
    \end{array}\right.\nonumber \\
   b_{jk}^{-} &=& \left\{\begin{array}{ll}     
     L_{\nu_{jk}}(Z_{j+1,k}-Z_{j,k}),&j\mbox{ even},\\ 
     L_{\nu_{jk}}(Z_{j+1,k+1}-Z_{j,k} ),&j\mbox{ odd},
    \end{array}\right. \label{eq:abc-} \\
   c_{jk}^{-} &=& \left\{\begin{array}{ll}     
     L_{\nu_{jk}}(Z_{j,k}-Z_{j+1,k-1}). &j\mbox{ even},\\ 
     L_{\nu_{jk}}(Z_{j,k}-Z_{j+1,k} ),&j\mbox{ odd}.
    \end{array}\right.\nonumber
\end{eqnarray}
However, while these equations stand as written for values of $k$
giving triangles in the ``interior'' of our mesh, at the upper and
lower parts of the mesh we must take the $2\pi i$-periodicity into
account in order to conserve our requirement that $0\le l\le N-1$ in
every appearance of $W_{jl}$. The exceptions to these equations occur
when $k+1=N+1$ in (\ref{eq:triangles+}) and when $k-1=-1$ in
(\ref{eq:triangles-}).  When $k=0$ and $j$ is even we should use
$W_{j\pm1,N-1}-2\pi i$ in place of $ W_{j\pm1,-1}$, while when $k=N-1$ and
$j$ is odd we should write $W_{j\pm1,0}+2\pi i$ instead of
$W_{j\pm1,N}$. Referring to (\ref{eq:triright}) and (\ref{eq:trileft})
one sees that the exceptional equations are defined by
\begin{eqnarray} \label{eq:triangles++}
  a_{j0}^{\pm}W_{j0} + b_{j0}^{\pm}W_{j-1,0} + c_{j0}^{\pm}W_{j-1,1} &=& 
     -2\pi i c_{j0}^{\pm},\quad j \mbox{ even},  \nonumber\\
 a_{j0}^{\pm}W_{j0} + b_{j0}^{\pm}W_{j-1,0} + c_{j0}^{\pm}W_{j-1,1} &=& 
     2\pi i b_{j0}^{\pm},\quad j \mbox{ odd} . 
\end{eqnarray} 

The next step is to ``reflect'' these equations to the right
half-plane via (\ref{eq:defrho}).  We want $F$ to be symmetric in the
imaginary axis, $F=\varrho F\varrho$. Since the Beltrami derivative of
this composition is $\nu=\conj{\nu}\circ\varrho$, the prescription
(\ref{eq:nusymm}) indeed reflects $\nu$ appropriately to the right
half-plane as a discrete Beltrami differential. In other words,
consider a triangle $\tau^\pm_{jk}=\varrho(\tau^\pm_{-j,k})$ in the right
half $Z$-plane ($j\ge0$).  The image $F(\tau^{\pm}_{j,k})$ in the
$W$-plane, defined by $W_{jk}$ and two adjacent vertices, must be the
same as $\varrho(F(\tau^\pm_{-j,k}))$. It is easily seen that the
correspondence $\varrho(\tau^{\pm}_{jk})\to \varrho(F(\tau^{\pm}_{jk}))$
translates into equations of the same form as (\ref{eq:triangles+}),
(\ref{eq:triangles-}) with $j\ge0$ and with the coefficients
\begin{equation} \label{eq:abcconj}
  a_{jk}^\pm = \conj{a_{-j,k}^\pm}, \ 
  b_{jk}^\pm = \conj{b_{-j,k}^\pm}, \ 
  c_{jk}^\pm = \conj{c_{-j,k}^\pm}.
\end{equation} 

So far we have described how the triangles of the extended $Z$-triangulation
provide $4MN$ linear equations: $2MN$ from the left half-plane and another
$2MN$ via (\ref{eq:abcconj}).

\subsection{Boundary equations}\label{subsec:bdryeq}

Next we look at the right and left boundary conditions. We return for
a moment to the $z$-disk.  The smallest polygon of the mesh is formed
of points on the circle of radius $r_{-M}=\exp(R_{-M})$. Following the discussion
in \ref{subsec:bdrycond}, let $e_k$ be the images of these points
under the real-linear mapping $L_{\mu_0}$,
\[   e_k = L_{\mu_0}(r_{-M}e^{2\pi i k/N}) = r_{-M} L_{\mu_0}(e^{2\pi i k/N}) 
 ,\quad 0\le k\le N-1,
\]
where $\mu_0$ denotes the average value of $\mu(z)$ inside this circle.
Thus the $e_k$ lie on a small ellipse. Define
\begin{equation} \label{eq:defEk}
  E_k = \log e_k
\end{equation}
with $0\le \arg E_k<2\pi$.  It might seem natural to use the values
$E_k$, $\varrho(E_k)$ as ``boundary values'' simply by adding
equations $W_{-M,k}=E_k$, $W_{M,k}=-\conj{E_k}$.  However, this will
create the difficulty which was mentioned earlier related to the
unknown conformal module of the region between the small ellipse and
its inversion, because the value of $\log r_{-M}$ in our construction,
which is the real part $R_{-M}$ of $E_k$, is arbitrary as far as
$\mu$ is concerned.  Instead of
that we want a condition which says that the image of the
circle of radius $r_{-M}$ is an unknown (complex nonzero) multiple of
the ellipse $\{e_k\}$.  In logarithmic coordinates, the condition is that the
image of the curve $\{Z_{-M,k}\}_k$ is a translate of the curve
$\{E_k\}_k$ by a complex constant.  Similar considerations apply
to the reflected curve $\{\varrho(E_k)\}_k$. The boundary equations which
achieve this are the $2(N-1)$ equations
\begin{eqnarray} \label{eq:bdry}
   W_{-M,k} - W_{-M,k-1}  &=& D_k, \nonumber \\ 
    W_{M,k} - W_{M,k-1}  &=&  \conj{D_k}   ,
\end{eqnarray}
where $D_k=E_k-E_{k-1}$ and $1\le k\le N-1$.  Note that the magnitude of
$r_{-M}$ does not influence the value of $D_k$.
 
Finally, for normalization of the solution we add one more equation,
\begin{equation}\label{eq:normalize}
 W_{0,0} = 0,
\end{equation}
which is self-symmetric. This says that $F(0)=0$, or equivalently,
$f(1)=1$.

\section{Statement and proof of main theorem}\label{sec:main}

The linear system outlined in section \ref{sec:DBE} has more
equations than variables: there are $n_{\rm v}=(2M+1)N$ unknowns
$W_{jk}$, $-M\le j\le M$, $0\le k\le N-1$, and $n_{\rm
  e}=4MN+2(N-1)+1$ equations. Therefore we will use the standard
Least-Squares approximation \cite{Bj} method to find a solution.  To
describe the system it is convenient to rename the variables 
in a single vector $V$ with
\begin{equation} \label{eq:VpWjk}
   V_p = W_{jk} 
\end{equation}  
where $p=p(j,k)$ is given by an arbitrary but fixed bijective
correspondence from the set of index pairs $\{(j,k)\}$ to the range
$1\le p\le n_{\rm v}$. We will write $V=\underline{W}$ or
$W=\underline{V}$ to indicate this renaming of the indices.  The
linear system now takes the form $A\underline{W}=B$ or
\begin{equation} \label{eq:AVB}
     A V=  B     
\end{equation}
where $A=(A_{np})$ is a complex $n_{\rm e}\times n_{\rm v}$ matrix and
$B=(B_n)$ is a complex vector of length $n_{\rm e}$.  When considering
the mesh $\{Z_{jk}\}$ as fixed, we will say that $(A,B)$ is the
\emph{associated linear system} to the collection of $\nu$-values $\{\nu_{jk}\}$
(recall that the coefficients depend both on $\nu_{jk}$ and $Z_{jk}$).
 
\subsection{Statement of theorem}
 
Our algorithm may be summarized briefly as follows:
\begin{enumerate}
\item Given a proposed Beltrami derivative $\mu$ in $\D$, choose the
  dimensions $M,N$ for a triangular mesh $\{Z_{jk}\}$ in the $Z$-plane
  and calculate the averages $\nu_{jk}$ of the pullback of $\mu(z)$ to
  the $Z$-plane via (\ref{eq:pullback}).
\item Calculate the coefficients of the linear system $(A,B)$
  associated to $\{\nu_{jk}\}$ as prescribed by equations
  (\ref{eq:triangles+}), (\ref{eq:triangles-}),
  (\ref{eq:triangles++}), (\ref{eq:bdry}), and (\ref{eq:normalize}).
\item Apply the method of Least Squares to calculate the approximation
  $V$ of the solution of the system $AV=B$, and arrange the entries of
  $V$ to form the mesh $\{W_{jk}\}=\underline{V}$.
\item Calculate $w_{jk} = \exp W_{jk}$ for $-M\le j\le 0$ and $0\le
  k\le N-1$.  The desired mapping is the picewise linear simplicial mapping such that $z_{jk} \mapsto w_{jk}$ where $z_{jk}=\exp Z_{jk}$.
\end{enumerate}
  
\begin{figure}[b !]
  \centering
  \pic{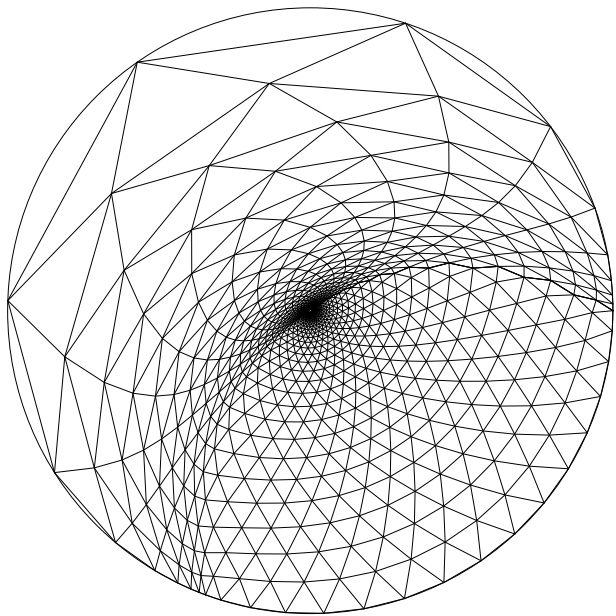}{1.5,0}{7.5}{}
  \caption{Different constants in upper and lower
    $w$-half-planes. Note the normalizations $f(0)=0$, $f(1)-1$.}
  \label{fig:splitdiskw}
\end{figure}
 
  \begin{figure}[th!]
  \centering
 \pic{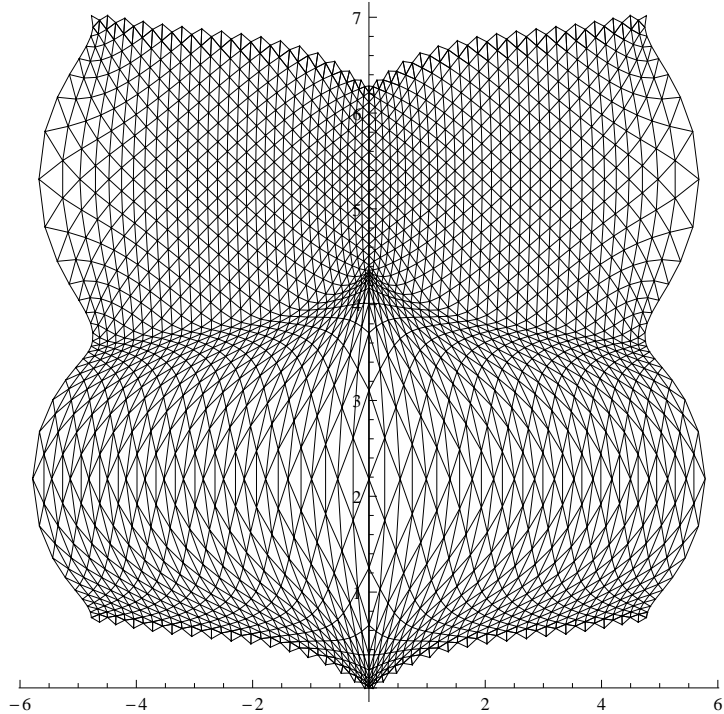}{1,-0}{8.5}{}
   \caption{Different constants in upper and lower $W$-half-planes.
Note the lifted ellipses at left and right extremes of the boundary.}
  \label{fig:splitdiskW}
\end{figure}
 
As an illustration we show in Figure \ref{fig:splitdiskw} the
$w$-triangulation for the Beltrami derivative defined by $\mu(z)=0.5$
when $\im z>0$, and $\mu(z)=0$ when $\im z\le0$.  This was calculated
with $(M,N)=(32,32)$. Therefore there are many very small triangles
which cannot be seen in the picture, in particular those adjacent to
the small bounding ellipse.  Note how the image triangles from the
lower half-plane appear to be equilateral, and there is a clear
dividing line between them and the stretched triangles from the upper
half-plane. The stretched triangles are not mutually similar, even
though $\mu$ is constant, since the similarity class of an image
triangle depends upon both the value of $\mu$ and the slope of the
base of the domain triangle.  In Figure \ref{fig:splitdiskW} we show
the logarithmic $W$-domain. In this picture the equilateral triangles
are in the upper part, with imaginary parts approximately from $\pi$
to $2\pi$ (near the imaginary axis, they occupy a smaller range).
 
Our main result is formulated as follows. Recall that the real parts
$R_j$ of the vertices of the logarithmic meshes are given by
(\ref{eq:Rj}).  

\begin{theo}\label{th:main}
  Let $\mu$ be a $C^1$ function in $\D$ with $\|\mu\|_\infty<1$.
  Let $M_s,N_s\to\infty$ as $s\to\infty$, where these sequences
  satisfy
\begin{equation}   \label{eq:MNgrowth}
   c_1 N_s\log N_s \le  M_s \le c_2N_s \log N_s
\end{equation}
for constants $c_1,c_2$ where $c_1>1/(\pi\sqrt{3})$. Then
\begin{itemize}
\item[i.] For large $s$, the points $\{z_{jk}^{(s)}\}$ and the points
  $\{w_{jk}^{(s)}\}$ produced by the algorithm form the vertex sets of
  isomorphic triangulations $\T_z^{(s)}$ and\/ $\T_w^{(s)}$ of the
  unit disk\/ $\D$. Further, any fixed compact set
  $K\subset\mbox{int}\;\D$ is contained in the supports of\/
  $\T_z^{(s)}$ and of\/ $\T_w^{(s)}$ for large $s$.
\item [ii.] Let $f^{(s)}$ denote the piecewise-linear mapping of
  $\T_z^{(s)}$ to $\T_w^{(s)}$ which sends $z_{jk}^{(s)}$ to
  $w_{jk}^{(s)}$ as given by Corollary \ref{coro:implicitmu}.  Then
  the mappings $f^{(s)}$ converge to the solution $f$ of the Beltrami
  equation (\ref{eq:beltrami}), normalized by $f(0)=0$, $f(1)=1$,
  as $s\to\infty$, uniformly on compact subsets of\/ $\D$.
\end{itemize}
\end{theo}

\subsection{Lemmas}

In the Lemmas \ref{lemm:AX0} to \ref{lemm:sym}, $(A,B)$ is the
associated linear system for some fixed collection $\{\nu_{jk}\}$ 
satisfying $|\nu_{jk}|<1$, with mesh dimensions $M$ and $N$.

Since there are many more equations than variables involved, it is not
surprising that the system is overdetermined:
\begin{lemm} \label{lemm:AX0}
 If $X\in\C^{n_{\rm v}}$ and $AX=0$, then $X=0$.
\end{lemm}
\proof We examine the role of different rows of $A$ in $AX=0$.
Recall that we write $\underline{X}_{jk}$ for $X_{p(j,k)}$.  First
note that there appear equations (\ref{eq:bdry}) for the ``lifted
ellipse,'' with zero on the right-hand side in place of $D_k$ and
$\conj{D_k}$. These equations say that $\underline{X}_{-M,k+1}-\underline{X}_{-M,k}=0$; i.e., all of the
$\underline{X}_{-M,k}$ are now of the same value, say $c$. The equations
(\ref{eq:triangles+}), (\ref{eq:triangles++}) for the
rightward-pointing triangles with $j=-M+1$ give us
\[ a_{-M+1,k}^{+}\underline{X}_{-M+1,k} + b_{-M+1,k}^{+}\underline{X}_{-M,k} +
   c_{-M+1,k}^{+}\underline{X}_{-M,k'} = 0
\] for suitable $k'$, and using what we have just proved for $j=-M$
and the facts that $a_{-M+1,k}^{+}\not=0$ (recalling the remark after
Proposition \ref{prop:Bdef}) and $ a_{-M+1,k}^{+} + b_{-M+1,k}^{+} +
c_{-M+1,k}^{+} = 0 $ we deduce that $\underline{X}_{-M+1,k}=c$ for all $k$.
Continuing this way we have $\underline{X}_{-j,k}=c$ for all $k$ and $-M\le j\le
0$. The normalization equation $\underline{X}_{0,0}=0$ says that $c=0$ and then
the symmetry gives $\underline{X}_{jk}=0$ for all $j,k$. Thus $X=0$.\qed

\begin{lemm} \label{lemm:AXn0}
  Let $\{X^{(n)}\}$ be such that $AX^{(n)}\to0$.  Then $X^{(n)}\to0$.
\end{lemm}
\proof
Consider a subsequence of $\{X^{(n)}\}$ which  converges to a limit
$X$. By continuity $AX=0$.  By Lemma \ref{lemm:AX0}, $X=0$.  Thus we
see that every convergent subsequence of  $\{X^{(n)}\}$ converges to 0.
If $\{X^{(n)}\}$ is bounded, then it indeed has a convergent subsequence,
so it follows that $X^{(n)}\to0$ as desired.

Suppose then that $\{X^{(n)}\}$ is not bounded.  The maximum absolute value
$|X^{(n)}_{p}|$ of an entry of  $X^{(n)}$  occurs infinitely often for some
fixed index $p_0$.  On the corresponding subsequence we have
$|X^{(n)}_{p_0}|\to\infty$. Let
\[ Y^{(n)} = \frac{1}{X^{(n)}_{p_0}} X^{(n)},
\]
so $|Y^{(n)}|=1$. Also
\[ AY^{(n)} = \frac{1}{X^{(n)}_{p_0}} AX^{(n)} \to 0
\]
on the subsequence, where $|X^{(n)}_{p_0}|>1$ large $n$.  Since
$\{Y^{(n)}\}$ is bounded, by the previous paragraph $Y^{(n)}\to0$,
which contradicts $|Y^{(n)}|=1$.  Therefore this case does not
occur. \qed

\begin{lemm} \label{lemm:sym} For any $W=\{W_{jk}\}$ ($-M\le j\le M$,
  $0\le k\le N-1$), the symmetry relation
\[ A\underline{\rev{W}} = A\underline{\varrho(W)}
\] 
holds, where $\rev{W}_{jk}=W_{-j,k}$,
$\varrho(W)=\{\varrho(W_{jk})\}$ and $\rho$ is defined by
$(\ref{eq:defrho})$.  Let $V$ be the solution of the linear system
(\ref{eq:AVB}) produced by the method of Least Squares. Then the
entries of $W=\underline{V}$ satisfy the symmetry
$W_{-j,k}=\rho(W_{jk})$. In particular, the central values $W_{0k}$ are purely
imaginary.
\end{lemm}

\proof The relation $A\underline{\rev{W}} = A\underline{\varrho(W)}$
follows immediately from the previous symmetry relations such as
(\ref{eq:abcconj}). The method of Least Squares produces $V$ which
minimizes the $L_2$-norm of the residual $\|AV-B\|_2$.  By Lemma
\ref{lemm:AX0}, $A$ has full column rank, so this optimal $V$ is
unique \cite{Bj}.  Therefore
$\rev{W}=\varrho(W)$. \qed

In the following, recall the distance between similarity classes
of triangles discussed in \ref{subsec:similarity}.
  
\begin{lemm} \label{lemm:triangulation} Let
  $f\colon\overline{\D}\to\overline{\D}$ be a $C^1$-diffeomorphism of
  the closed disk (i.e., the restriction of a diffeomorphism of larger
  domains). Let $0<c<1$. Then every sufficiently fine triangulation
  $\T_z$ of $\D$, formed of triangles which are within $c$ of being
  equilateral, has the property that the images under $f$ of its
  vertices form the vertices of an isomorphic triangulation $\T_w$.
\end{lemm}

\proof We may suppose without loss of generality that $f$ preserves
orientation. It is sufficient to show that all triangles of $\T_z$ are
sent to triangles of the same orientation. Indeed, this implies that
the induced PL-mapping $\hat f$ preserves the triangle structure
locally, and hence is a local homeomorphism. Further, it is easily
seen that $\hat f\colon|\T_z|\to|\T_w|$ satisfies the path-lifting
property, and thus is a covering map of simply connected regions and
hence a homeomorphism, so the triangulations are isomorphic.
 
If the affirmation were false, there would be a sequence
$\{\T^{(n)}\}$ of triangulations of $\D$ such that all triangles
$\tau\in\T^{(n)}$ have diameter no greater than $1/n$, and for every
$n$ there is some $\tau_n\in\T^{(n)}$ such that $f$ reverses the
orientation of the vertices of $\tau_n$.  On a subsequence there is a
limit point $\tau_n\to z_0\in\D$. By the restriction on the similarity
classes, there is a uniform upper bound $c'$ to the ratio of any two
sides of $\tau_n$. In a suitable neighborhood $U$ of $z_0$ the linear
approximation
\[  |f(z) - f(z_0) - (J_f|_{z_0})(z-z_0)| \le \frac{1}{3c'}|z-z_0|,
\]
is valid.  For large $n$, $\tau_n\subseteq U$, and the inequality
implies that that $f$ sends the vertices of $\tau$ near to their
images under the affine-linear mapping $f(z_0)+(J_f|_{z_0})(z-z_0)$,
contradicting the property of reversing the orientation. This proves
the affirmation.  \qed

\subsection{Proof of main theorem}
\newcounter{step}\def\step{\addtocounter{step}{1}\noindent{\bf(\arabic{step}.)} }

We divide the proof into several steps.

\step By standard approximation arguments we may replace $\mu(z)$ with
$\mu(rz)$ for $r<1$ arbitrarily close to 1, and thus we may suppose
that $\mu$ is $C^1$-smooth in a neighborhood of $\overline{\D}$.  By
Lemma \ref{lemm:triangulation}, the solution $f=f_\mu$ of $(\partial
f/\partial\overline{Z})/(\partial f/\partial Z)=\mu$, which is indeed
of class $C^2$, sends sufficiently fine triangulations to
triangulations.  Statement (i) of the theorem will follow when we
prove that the algorithm produces a sufficiently good approximation to
$f$.   

The pullback $\nu$ of $\mu$ given by (\ref{eq:pullback}) is $C^1$ in
the closed left half-plane and is periodic of period $2\pi i$, while
the solution $F=F_v$ of $(\partial F/\partial\overline{Z})/ (\partial
F/\partial Z)=\nu$ is $C^2$ there, with the limiting conditions
$F(-\infty)=-\infty$, $F(0)=0$.  Note that in general the extension of
$F$ by reflection in the imaginary axis is not $C^1$ on the axis.

The growth condition (\ref{eq:MNgrowth}) implies that
\begin{equation}\label{eq:rshrink}
  r_{-M} < \frac{1}{N},
\end{equation}
so further 
\begin{equation} \label{eq:zshrink}
   |z_{-M,k}-z_{-M,k-1}| < \frac{2\pi}{N}r_{-M}  = O\left(\frac{1}{N^2}\right)
\end{equation}
where $M=M_s$, $N=N_s$, $s\to\infty$.
 
\step
We recall that the statement ``$AV=B$'' in step 3 of the Algorithm
translates roughly into the statement that ``the PL-mapping $Z\mapsto
W$ is $\{\nu_{jk}\}$-conformal''. However, the algorithm only produces
the Least-Squares approximation for $(A,B)$, which we will call $V'$;
i.e., the $L_2$-norm $\|R'\|_2$ of the residual vector
\begin{equation} \label{eq:R'}
 R'=AV'-B 
\end{equation}
is the smallest possible.

We now restrict our attention to the mesh $\T_{MN}$. Consider the 
vector $V$ defined by
\begin{equation}
 \underline{V}= W=F(\T_{MN}). 
\end{equation}
which contains the images of the vertices under the true $\nu$-conformal
mapping $F$.  By Lemma \ref{lemm:triangulation}, $W$ is a
triangulation when the values $M,N$ are large enough.  Fix such $M,N$,
and let $\nu_{MN}=\{(\nu_{MN})^\pm_{jk}\}$ denote the collection of
 average values of the function $\nu$ on the triangles of
$\T_{MN}$.  Recall that the
associated linear system $(A,B)=(A_{MN},B_{MN})$ used in the algorithm
is defined in terms of the values of $\nu_{MN}$. 
 
Let $F^*_{MN}$ denote the PL-mapping on the support of $\T_{MN}$
defined by the condition $\T_{MN}\to \underline{V}$. Thus by
construction $F^*_{MN}$ coincides with $F$ on the vertices of
$\T_{MN}$ (both map to $W$), but the Beltrami derivative of $F^*_{MN}$
is constant on each triangle. We will write $\nu^*_{MN}$ for this
collection of constants.  
 
Let $(A^*,B^*)$ be the associated linear system  to the discrete Beltrami
coefficient $\nu^*_{MN}$. We consider the following for
fixed values of $M,N$:
\begin{eqnarray} 
  A V' - B\ \  &=& R'. \label{eq:AV'BR'} \\
  A^* V  - B^* &=& \varepsilon.   \label{eq:A*V*B*}  
\end{eqnarray}

\step The vector $ \varepsilon=\varepsilon_{MN}$ defined in
(\ref{eq:A*V*B*}) is described as follows. The entries of
$\varepsilon$ corresponding to the the triangle equations are 0
because $F^*_{MN}$ is $\nu^*_{MN}$-conformal and takes $\T_{MN}$ to
$\underline{V}$. Therefore the only nonzero values in $\varepsilon$
are in the positions corresponding to the boundary equations.  By
(\ref{eq:bdry}) the values in $A^*V$ are $F(Z_{-M,k})-F(Z_{-M,k-1})$
while the values in $B^*$ are $E_k-E_{k-1}$, so in the corresponding
positions $\varepsilon$ contains the value
\[  (F(Z_{-M,k})-E_k) - (F(Z_{-M,k-1})-E_{k-1}).
\]
To estimate this we recall $f$ is of class $C^2$ and use the
approximation
\[  f(z) = (J_f|_ 0)(z) + O(|z|^2)
\]
at the origin, where the Jacobian is given by $J_f|_0=H_{a,0}\circ L_{\mu(0)}$ for
some complex $a\not=0$.  From (\ref{eq:zshrink}) we see that
\[ \frac{f(z_{-M,k})e_{k-1}}{f(z_{-M,k-1})e_{k}} = 
   \frac{e_ke_{k-1}+O(r_{-M}^3)}{e_ke_{k-1}+O(r_{-M}^3)} = 1 + O(r_{-M}) 
\]
Taking logarithms we conclude that  the nonzero entries of $\varepsilon$
shrink at least as fast as $O(r_{-M})$. There are at most
$2N$ nonzero entries, so (\ref{eq:rshrink}) gives
\begin{equation}\label{eq:varepsilon}
 \|\varepsilon_{MN}\|_2 = O\left(\left((2N)\frac{1}{N^{2}}\right)^{(1/2)}\right) 
 \to 0.
\end{equation}

\step In practice, one applies the algorithm by using average values
\begin{equation}   \label{eq:nuMN}
  \nu_{MN}(T) = \frac{\nu(Z_1)+\nu(Z_2)+\nu(Z_3)}{3}.
\end{equation}
for $T=(Z_1,Z_2,Z_3)\in\T_{MN}$. It is a simple exercise to verify that this differs from
the value $(\int_T\nu\,dx\,dy)/\mbox{area}(T)$ by an amount which tends to zero as $O(\delta^2)$, 
 where
$\delta=\mbox{diam } T$. Since $\nu$ is differentiable,
$\nu(Z_i)=\nu(Z_0)+(J_\nu|_{Z_0})(Z_i-Z_0)+ O(\delta^2)$, $i=1,2,3$,
where we take $Z_0=(Z_1+Z_2+Z_3)/3$.  It is immediate from
(\ref{eq:nuMN}) that
\[ \nu_{MN}(T) = \nu(Z_0)+O(\delta^2)
\]
as the triangulation is refined and $Z_0$ always refers to the center of the triangle $T$.

On the other hand, consider
\begin{eqnarray} 
  W_i &=& F(Z_i) \ = \ F^*(Z_i) \ = \ W_0 + (J_f|_{Z_0})(Z_i-Z_0) + O(\delta^2) \nonumber\\
  &=&   W_0 + H_{a,b}\circ L_{\nu(Z_0)}(Z_i-Z_0) + O(\delta^2)  \label{eq:jacobian}
\end{eqnarray}
(where the constants $a,b$ depend on $Z_0$) 
and recall that $\nu^*_{MN}(T)$ is given by Corollary \ref{coro:implicitmu}, i.e.,
\begin{eqnarray*}
 \nu^*_{MN}(T) &=& \frac{}{}
  -\frac{(Z_2-Z_1)(W_3-W_1) - (Z_3-Z_1)(W_2-W_1)}
      {(\conj{Z_2}-\conj{Z_1})(W_3-W_1) - 
       (\conj{Z_3}-\conj{Z_1})(W_2-W_1)} \\
&=& \nu(Z_0) + O(\delta)
\end{eqnarray*}   
as is seen by applying (\ref{eq:defL})--(\ref{eq:defH}) in
(\ref{eq:jacobian}) and cancelling.  We conclude that
\begin{equation} \label{eq:nu*nuMN}
  |\nu^*_{MN} - \nu_{MN}| =   O(\delta) = O\left(\frac{1}{N}\right).
\end{equation}

\begin{figure}[t!]
 \pic{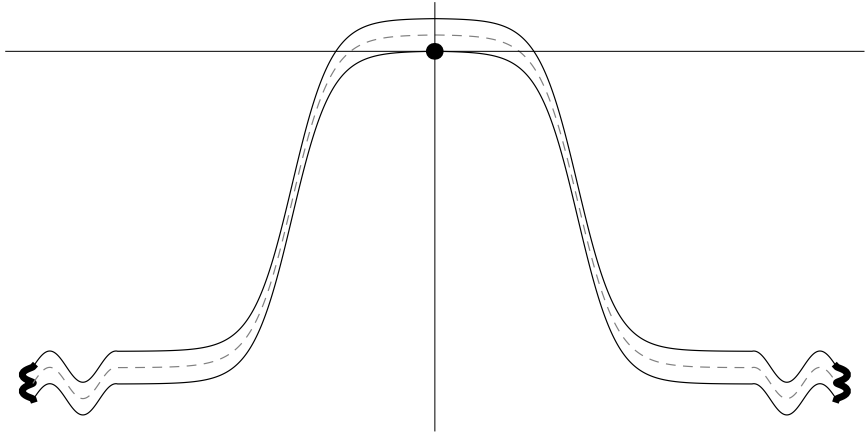}{5.4,0}{5}{scale=.9}
 \caption{$W$-image containing points far from the
origin.  The lifted ellipse and its reflection are drawn as thicker curves
for emphasis.
} \label{fig:bigimage}
\end{figure}

\step  
The next step is to verify that 
\begin{equation}  \label{eq:Vbound}
  \|V\|_\infty = O(\log N).
\end{equation}
The support $|\T_{MN}|$ of the $Z$-triangulation is a rectangle of
width $R_{M}\approx\log N$ and height $2\pi$. By (\ref{eq:A*V*B*}),
(\ref{eq:varepsilon}) we see that $\|A^*V-B^*\|_\infty$ can be made
arbitrarily small by refining the mesh. This says that the triangles
of $\underline{V}$ are close to being $\nu^*_{MN}$-conformal, so in
particular the PL-mapping $F^*_{MN}\colon T_{MN}\to \underline{V}$ is
now known to be a local homeomorphism.  We may think of the image as a
Riemann surface extended over a region of the $W$-plane (Figure
\ref{fig:bigimage}).  This image is a possibly non-schlicht
topological quadrilateral whose ``vertical'' sides are the lifting of
the small ellipse together with its reflection in the imaginary
axis. The ``horizontal'' sides are two curves shifted from one another
by approximately $2\pi i$.  If (\ref{eq:Vbound}) were false, this
quadrilateral would contain points $W_s$ such that
$\rho_s:=|W_s|/\log N_s\to\infty$. This would imply that all curves joining
the two vertical sides (see dotted curve in Figure \ref{fig:bigimage})
have Euclidean length at least $\rho_s\log N_s$.  It follows from
\cite[Lemma 4.1]{LV}, that the conformal module of this quadrilateral
must also grow at least as fast as $\rho_s\log N_s$. Since $|\T_{MN}|$ has
conformal module $O(\log N_s)$, the quasiconformal mappings $F^*_{MN}$ must
have arbitrarily large maximal dilatation; i.e., their Beltrami
derivatives must have absolute value near to $1$ at some point. This
contradicts the fact that the Beltrami derivative is $\nu^*_{MN}$,
which is bounded away from $1$ since it is close to $\nu$.  Therefore
(\ref{eq:Vbound}) holds as claimed.
 
\step
Now  write
\[  |A^*-A| = \sup |(A^*)_{jk}-(A)_{jk}|, \ |B^*-B| = \sup |(B^*)_k-(B)_k|.
\]
We observe that entries of $(A,B)$ and $(A^*,B^*)$ in a given position
come from identical explicit formulas, with data $\nu^*_{MN}$ and
$\nu_{MN}$ respectively. By Corollary \ref{coro:eqntri}, entries in
$A$, $A^*$ resulting from triangle equations are of the form
$L_\nu(Z_i-Z_j)$, $L_{\nu^*}(Z_i-Z_j)$, respectively, where a given
triangle is referred to by $(Z_1,Z_2,Z_3)$ as in (\ref{eq:nuMN}), and
$\nu$, $\nu^*$ refer to the constant values assigned to that
particular triangle.  The difference is
\begin{eqnarray*}
 L_\nu(Z_i-Z_j)-L_{\nu^*}(Z_i-Z_j) &=&
  \frac{2(\nu^*-\nu)\im(Z_i-Z_j)}{(1+\nu)(1+\nu^*)}  \\
   &=& O\left(\frac{1}{N^2}\right)
\end{eqnarray*}
by (\ref{eq:nu*nuMN}), since $|\nu|$ is bounded away from 1.  The
entries of $A,A^*$ resulting from boundary equations do not depend on
$\nu$, $\nu^*$.  We arrive at
\begin{equation} \label{eq:ABA*B*}
  |A^*-A|=O(1/N^2), \quad|B^*-B|=O(1/N),
\end{equation}
the latter estimate resulting from a simple calculation based on
(\ref{eq:bdry}).
 
Next we observe that by (\ref{eq:A*V*B*}).
\begin{eqnarray*}
  \|AV-B\|_2 &\le& \|AV-A^*V\|_2 + \|A^*V-B^*\|_2 + \|B^*-B\|_2 
    \nonumber\\
  &=&\|(A-A^*)V\|_2 + \|\varepsilon\|_2 + \|B^*-B\|_2 \label{eq:AV-B} .
\end{eqnarray*}
By construction, each row of $A$ or $A^*$ contains at most three
nonzero entries, as only three variables appear in equations
(\ref{eq:triangles+}), (\ref{eq:triangles-}), etc. The entries of
$(A-A^*)V$ are of the form $\sum_p(A_{pq}-A^*_{pq})V_p$, where for
each $q$, at most three of the summands are nonzero. Thus
$\|(A-A^*)V\|_\infty\le3|A-A^*|\,\|V\|_\infty=O(\log N/N^2)$ by
(\ref{eq:Vbound}) and (\ref{eq:ABA*B*}).  This vector has ${\rm n_e}$
elements, so we may estimate its $L_2$-norm,
\begin{eqnarray*}
  \|(A-A^*)V\|_2 &=& O\left( \left({\rm n_e}(\frac{\log N}{N^2 })^2\right)^{1/2}
  \right) \ = \ O\left( \left(N^3\frac{1}{N^4}\right)^{1/2}\log N\right) \\
  &=& O(N^{-1/2}\log N) \to 0.
\end{eqnarray*}

A similar calculation shows that $\|B^*-B\|_2\to0$. 
By minimality of $\|AV'-B\|_2=\|R'\|_2$ (recall (\ref{eq:R'})) and by
(\ref{eq:varepsilon}), we have  
\[ \|R'\|_\infty<\|R'\|_2\le\|AV-B\|_2. \] 
 Therefore we have proved
that $\|R'\|_\infty\to0$. We will have no further need of the
$L_2$-norm in the discussion.
 
\step By  (\ref{eq:AV'BR'}) and  (\ref{eq:A*V*B*}),
\begin{eqnarray}
 \|A^*(V-V')\|_\infty  &\le& \|A^*V-B^*\|_\infty + \|B^*-B\|_\infty + \|B-AV'\|_\infty
 \nonumber \\
     &&\ \ +\ \|AV'-A^*V'\|_\infty \nonumber\\
   &=&    \|\varepsilon\|_\infty+\|B^*-B\|_\infty+\|R'\|_\infty  \nonumber \\
     &&\ \ +\  \|(A-A^*)V'\|_\infty. \label{eq:A*VV'}
\end{eqnarray}
The fact $\|AV'-B\|_\infty\to0$ tells us that the PL-mapping
$\T_{MN}\to W'=\underline{V'}$ has bounded dilatation. As we showed in
(\ref{eq:Vbound}) for $V$, it follows also that $\|V'\|_\infty =
O(\log N)$.  Thus we may conclude that $\|(A-A^*)V'\|_\infty\to0$
(again we need the fact that the rows of $A^*$ have no more than three
nonzero entries), so (\ref{eq:A*VV'}) implies that
$\|A^*(V-V')\|_\infty\to0$.  By Lemma \ref{lemm:AXn0} (applied to
$A^*$ in place of $A$), we see that $\|V-V'\|_\infty $ (as determined
by $M_s,N_s$) tends to zero as $s\to\infty$. This says that the points
$W'=\underline{V'}$ produced by the algorithm differ by an arbitrarily
small amount from the image vertices under the true $\nu$-conformal
mapping $W=F(Z)$.

\step Finally, we apply the exponential mapping via (\ref{eq:WFZ}) and
(\ref{eq:expZexpW}) to obtain the sequence of PL-mappings
$f^{(s)}\colon z\mapsto w$ of subdomains which exhaust the unit disk
$\D$, produced by the algorithm for the meshes determined by
$(M_s,N_s)$. Let $\epsilon>0$, and consider the annulus
$\D_\epsilon=\{\epsilon<|z|<1\}$. The image $f(\D_\epsilon)$ is
approximately the part of $\D$ outside of a small ellipse. By
construction, the correspondence $F\colon Z\mapsto W$ extends by $2\pi
i$-periodicity to a mapping of the left half-plane $\{\re Z\le0\}$ to
$\{\re W\le0\}$, as do all the approximants $F^{(s)}$.  The extended
quasiconformal mapping $F$ is uniformly continuous on the band
$\{-1/\epsilon\le \re Z\le0\}$, so the diameters of the $W$-triangles
in the images of triangles contained in this band tend uniformly to
zero as $s\to\infty$. Consider such a triangle $(W_1,W_2,W_3)$ and its
image $(w_1,w_2,w_3)$ where $w_i=\exp W_i$. We compare the angle at
$W_1$, which is $\arg(W_3-W_1)/(W_2-W_1)$, with the angle
\[ \arg\frac{w_3-w_1}{w_2-w_1}= \arg\left(
   \frac{e^{(W_1+W_3)/2}}{e^{(W_1+W_2)/2}}\cdot\frac{\sin(W_3-W_1)}{\sin(W_2-W_1)}\right)
\] 
at $w_1$.  The first factor on the right-hand side is
$e^{(W_3-W_2)/2}\to1$ since $|W_3-W_2|\to0$. The second factor tends
to $(W_3-W_1)/(W_2-W_1)$, again uniformly in the band.  As a
consequence, the exponential mapping sends $W$-triangles to
approximately similar $w$-triangles (thus in particular respecting the
orientation), and the procedure provides a triangulation of
$f^{(s)}(\D_\epsilon)$ for which the induced PL-mapping
$f^{(s)}|_{\D_\epsilon}$ is approximately $\mu$-conformal. The limit
as $s\to\infty$ is $\mu$-conformal, fixes $z=1$, and hence coincides
with the mapping of doubly connected domains
$f|_{\D_\epsilon}\colon\D_\epsilon\to f(\D_\epsilon)$. Since
$\epsilon$ is arbitrary, we conclude that $f^{(s)}\to f$. This
completes the proof. \qed

\section{\label{secnumres}Numerical Results\label{secnumerical}}
 
All of the calculations have been done with machine precision in
\emph{Mathematica} on a standard laptop computer of approximately 1GH.
We have not found any examples where more precision will make a
difference. The \emph{Mathematica} routine \texttt{LeastSquares}
handles sparse matrices \cite{Bj,Tr}, a data structure which registers
only the nonzero entries appearing in a matrix or vector.   
  
In the first several examples we compare the results produced by our
algorithm with an exact formula for the quasiconformal mapping under
consideration.

\example A simple test case is for constant $\mu(z)=c$.  There is an explicit
formula \cite{Sz} for the conformal mapping to $\D$ from an ellipse
with semimajor and semiminor axes of lengths $a$, $b$ ($a^2-b^2=1$) and foci
at $\pm1$,
\begin{equation} \label{eq:ellipticint}
  w=\sqrt{k}\, {\rm sn}\,(\frac{2K}{\pi} \sin^{-1} u ; k^2) 
\end{equation}
where the Jacobi elliptic function modulus $k$ is related to the
complete elliptic integral $K$ and the Jacobi theta functions by the
formulas
\[ q=(a+b)^{-4}=e^{-\pi K(1-m)/K(m)}, \]
\[ k = \sqrt{m} = \left( \frac{\theta_2}{\theta_3}\right)^2, \] with
notation from \cite{WW}. Note that the image of the circle $|z|=1$
under the mapping $L_\mu$ is an ellipse with semiaxes $1$,
$(1-|\mu|)/(1+|\mu|)$ slanted in the directions $(1/2)\arg \mu$,
$(1/2)(\arg\mu+\pi)$ respectively, modulo $\pi$. This ellipse is sent
by the conformal linear mapping $H_{1/(2\sqrt{\mu}),0}$ to the ellipse
with semiaxes $a,b$.  Then via (\ref{eq:ellipticint}) this is transformed
conformally to the unit disk.
 
\begin{figure}[!ht]
\pic{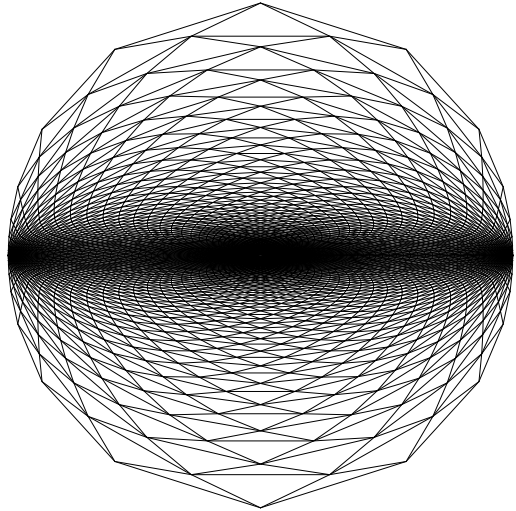}{8,0}{6}{} 
\caption{Image for constant Beltrami derivative $\mu=0.5$ and
  $(M,N)=(52,64)$. Observe the ``crowding phenomenon'' at the
  boundary.}
\label{fig:const}
\end{figure}

It is well known that a conformal mapping from an ellipse to a disk
tends to ``crowd'' boundary points near the images of the endpoints of
the major axis. The crowding, or maximum ratio of separation of $N$
points sent to the $N$th roots of unity, increases exponentially as a
function of the aspect ratio $a/b$ \cite[section 2.6]{DT}. In
contrast, the affine mapping $L_\mu$ which we combined with this
conformal mapping produces little crowding. The combined effect is a
great deal of crowding near $w=1$, as can be perceived in Figure
\ref{fig:const}.
 
The algorithm of Theorem \ref{th:main} was applied for the constant
Beltrami derivatives $\mu=0.1,\ 0.3,\ 0.5,\ 0.7$, and meshes defined
by $N=16$, $32$, $48$, $64$, $72$, $84$, with $M$ equal to the least
multiple of 4 no less than $N\log N/(\pi\sqrt{3})$, in view of
(\ref{eq:MNgrowth}). In the last case there are 24359 equations in
14196 variables. It took about 1.5 seconds to calculate the part of the
matrix in the left half-plane, and about 10 seconds to solve the full
set of equations.  Table \ref{tab:constants} presents the
maximum error over all $w_{kj}$ when these points are compared to the
images of $z_{kj}$ under the exact quasiconformal mapping described in
the preceding paragraph.  As is to be expected, the error increases
when the Beltrami derivative increases, but decreases when the mesh is
refined.  It was found that for $\mu=0.7$ and a rather coarse mesh
such as $(M,N)=(36,48)$, the image of $\T_{MN}$ is not a
triangulation, inasmuch as a few of the $w$-triangles near $\pm1$ are
improperly oriented. In spite of this fact, the values obtained for
the conformal mapping are not very far off.
 
\begin{figure}[ht!]
  \centering  
\pic{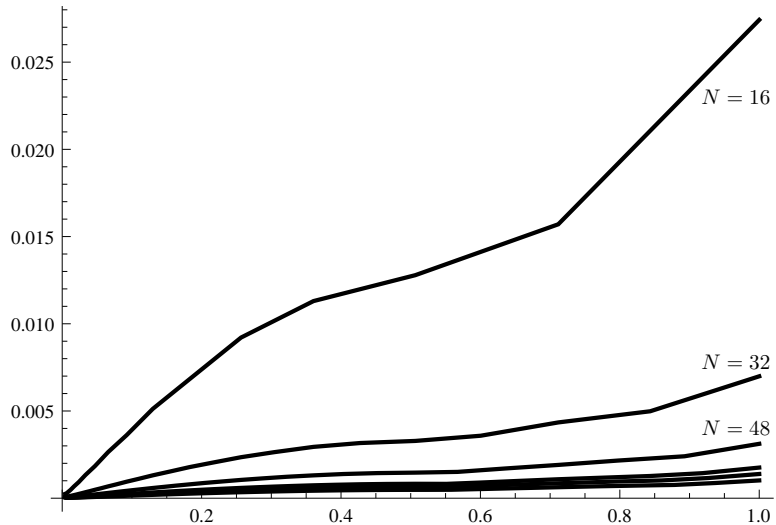}{0,0}{6.5}{scale=.8}
\caption{Numerical errors of algorithm for different
  values of $(M,N)$ with $\mu=0.3$.  The horizontal axis indicates the
  distance $r_j=|z_{jk}|$ of the $z$-points from the origin; the vertical axis
  gives the maximum discrepancy (over $k$) of the calculated value of
  $w_{jk}$ from the true value. }
  \label{fig:exampleaffine}
\end{figure} 
 
\begin{table}[!hb]
\begin{center}
\begin{tabular}{c||c|c|c|c|c|c}
$(M,N)$ & (12,16) & (24,32) & (36,48)& (52,64) & (60,72) & (72,84) \\ \hline
$\mu=$0.1 &0.012 & 0.0031 & 0.0014 & 0.0008 & 0.0006 & 0.0004\\
$\mu=$0.3 & 0.0274 & 0.007 & 0.0031 & 0.0018 & 0.0014 & 0.001\\
$\mu=$0.5 & 0.0615 & 0.0205 & 0.0109 & 0.0065 & 0.0051 & 0.0038\\
$\mu=$0.7 & 0.2439 & 0.1201 & 0.0856 & 0.0627 & 0.053 & 0.0412\\
\end{tabular}
\caption{The maximum of the absolute errors between the solutions and the real values of some constant Beltrami derivative and $M\approx N\log N/(\pi\sqrt{3})$. }
\label{tab:constants}
\end{center}
\end{table}

We give a further analisis of the variation of the error as a function
of the radius, for the particular value $\mu=0.3$.  Figure
\ref{fig:exampleaffine} shows the maximum error over $k$ in the
calculated value of $w_{jk}$ for each fixed $j$.  It is seen that the
error remains approximately constant for $r<0.7$ and then increases
rather sharply for $0.7<r<1$. Thus the maximum values in Table
\ref{tab:constants} are much higher than the average errors.  As noted
above, the maximum error, which occurs on the boundary, decreases as a
function of $N$.
    
\begin{figure}[!ht]
\pic{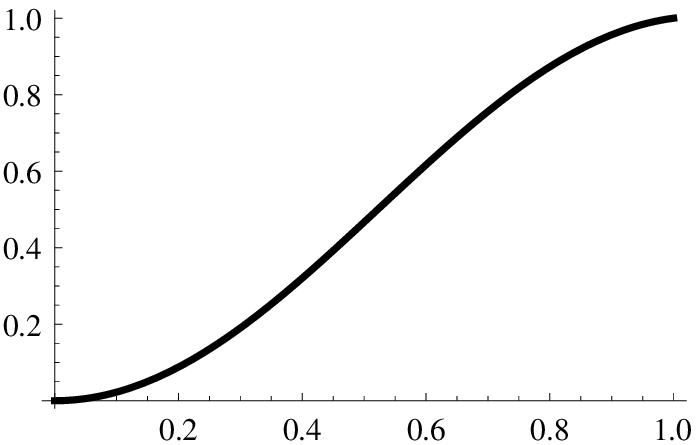}{5,0}{5}{scale=1}
\pic{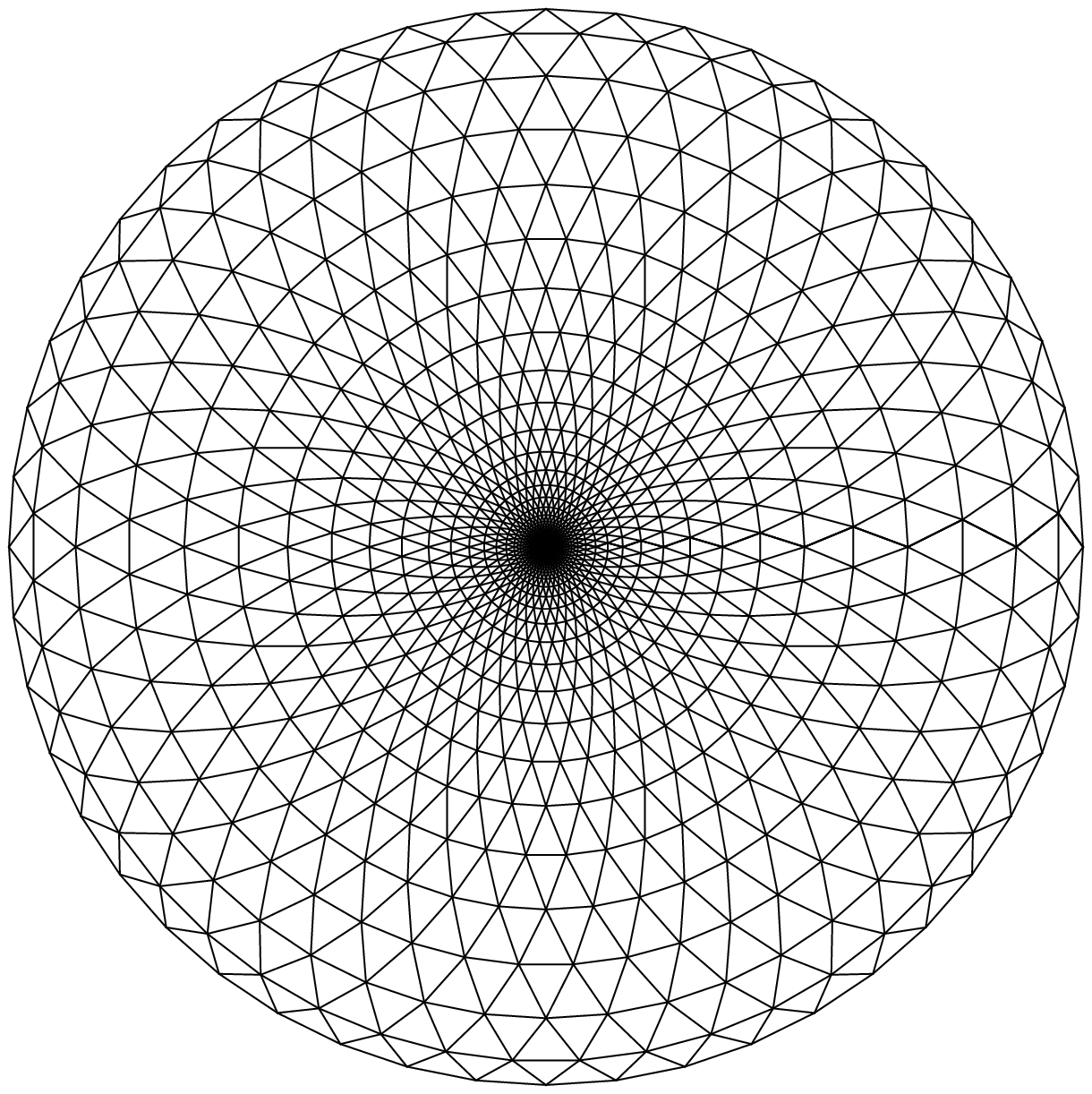}{11.9,0}{5}{scale=.40}
\caption{Radial function $\varphi$ of example 2 (left),
together with the induced rotationally symmetric image domain.}
\label{figradial} 
\end{figure}

\example \emph{Radial quasiconformal mapping.}  Let
$\varphi\colon[0,1]\to[0,1]$ be an increasing diffeomorphism of the
unit interval. Then the radially symmetric function
\begin{equation} \label{eq:fradial}
   f(z) = \varphi(|z|)e^{i \arg z} = \varphi(|z|)\frac{z}{|z|}
\end{equation}
has Beltrami derivative equal to
\begin{equation} \label{eq:muradial}
  \mu(z) = \frac{|z|\varphi'(z)/\varphi(z) - 1  }
                    {|z|\varphi'(z)/\varphi(z) + 1  } \frac{z}{\overline{z}}
\end{equation}
when $z\not=0$. As an illustration we will take 
\[ \varphi(r) = (1 - \cos 3r)/(1-\cos 3)
\]
as in Figure \ref{figradial}.  The resulting Beltrami derivative
satisfies $\|\mu\|_{\infty}=0.65$ approximately.

The domain points $z_{jk}$ on the real axis were selected, and the
values of $w_{jk}$ produced by the algorithm were compared with with
the true values $\varphi(|z_{jk}|)$. The results are given in Table
\ref{tab:radial}.  It was also observed that as in the previous
example, the errors increase as the radius increases.

\begin{table}[!htb]
\begin{center}
$\begin{array}{c||c|c|c|c|c|c}
$(M,N)$ & (12,16) & (24,32) & (36,48)& (52,64) & (60,72) & (72,84) \\ \hline
\rm Error &   0.0398 &  0.0135 &  0.0058  &  0.0034&   0.0027 &   0.0020 
\end{array}$
\caption{Maximum absolute error for radially symmetric with $\mu$
  mapping defined by (\ref{eq:muradial}). }
\label{tab:radial}
\end{center}
\end{table}

\begin{figure}[!ht]
\pic{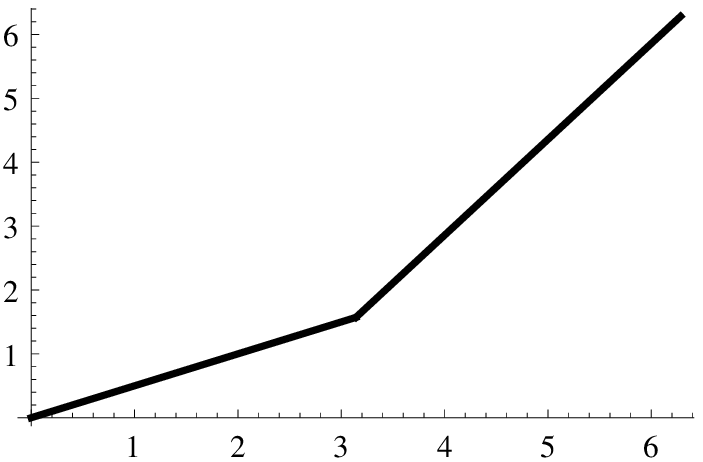}{5,0}{5}{scale=.9}
\pic{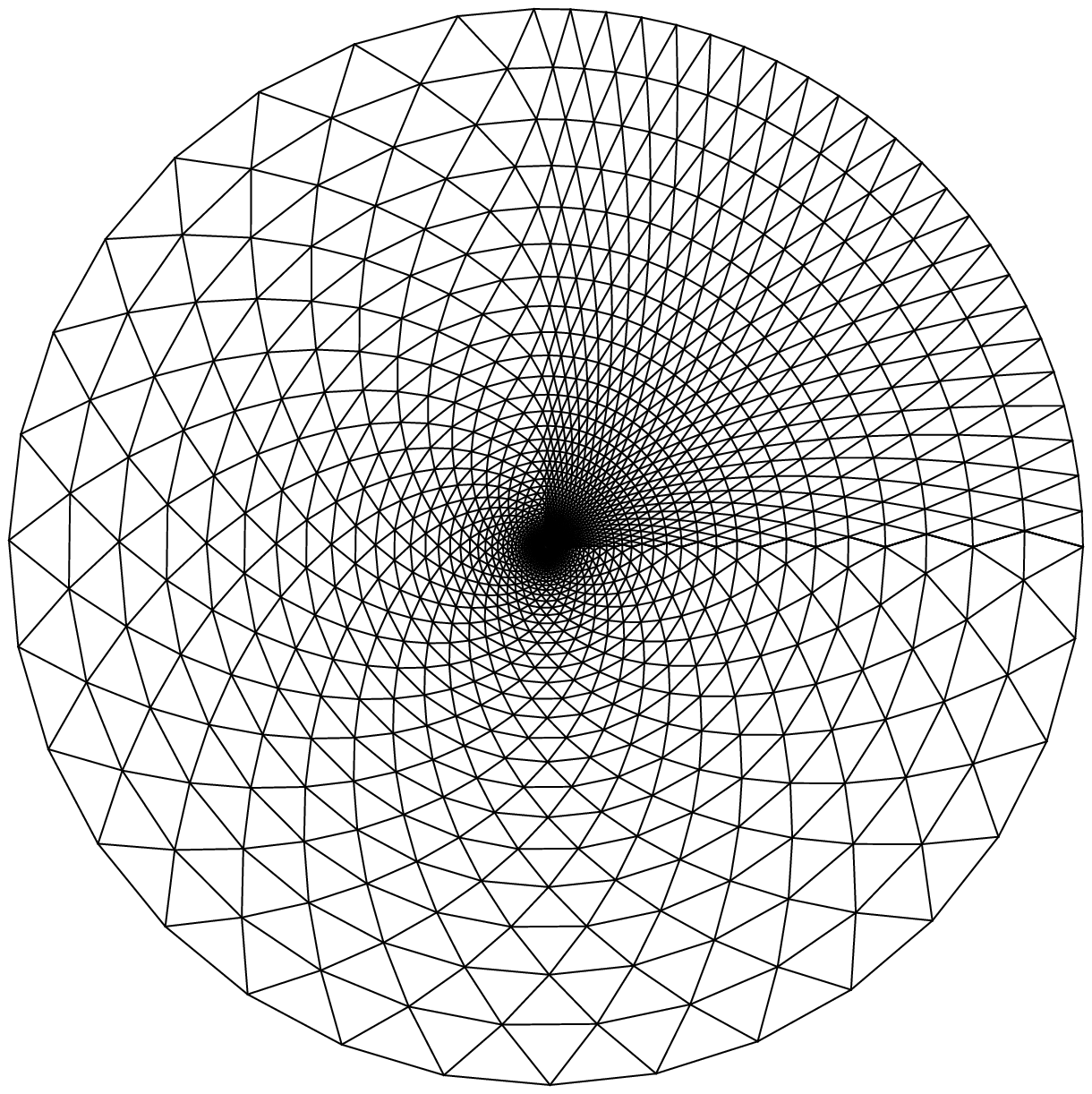}{11.9,0}{5}{scale=.40}
\caption{Angular function $\psi$ of example 3
(left), together with image domain under sectorial mapping (\ref{eq:fsectorial}).}
\label{fig:sectorial} 
\end{figure}
 
\example \emph{Sectorial quasiconformal mapping.}  In a similar
spirit, we let $\psi\colon[0,2\pi]\to[0,2\pi]$ be an increasing
diffeomorphism.  Write
$\widetilde\psi(e^{i\theta})=e^{i\psi(\theta)}$. Then the sectorially
symmetric function
\begin{equation} \label{eq:fsectorial}
   f(z) = |z|\,\widetilde\psi\left(\frac{z}{|z|} \right)
\end{equation}
has Beltrami derivative equal to
\begin{equation} \label{eq:musectorial}
  \mu(z) = \frac{1-\psi'(\theta)}{1+\psi'(\theta)}\,
           \frac{z}{\conj{z}}         
\end{equation}
when $z\not=0$. As an example we will take 
\[ \psi(\theta) =  \left\{ \begin{array}{ll}
         \frac{\theta}{2}, \quad& 0\le\theta\le\pi,\\
         \frac{\pi}{2}+\frac{3(\theta-\pi)}{2},& \pi\le\theta\le2\pi.
			 \end{array} \right. 
\]
as in Figure \ref{fig:sectorial}.  In this example $\mu$ does not 
satisfy the hypotheses of Theorem \ref{th:main} because it is not
continuous.  The arguments of the final boundary values on the
unit circle were compared with the true values $\psi(\theta)$; see
Table \ref{tab:sectorial}.

\begin{table}[!h]
\begin{center}
$\begin{array}{c||c|c|c|c|c|c}
$(M,N)$ & (12,16) & (24,32) & (36,48)& (52,64) & (60,72) & (72,84) \\ \hline 
\rm Error   & 0.0712 & 0.0362 & 0.0251 & 0.0193 & 0.0173 & 0.0150
\end{array}$
\caption{Maximum absolute error $|\psi(\theta)-f(e^{i\theta})|$ for sectorial mapping with $\mu$
    defined by (\ref{eq:musectorial}). }
\label{tab:sectorial}
\end{center}
\end{table}

\begin{figure}[b!]
\pic{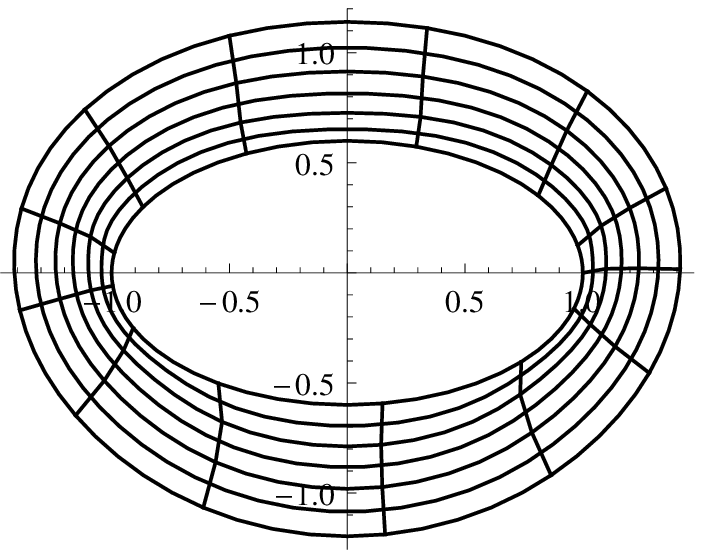}{5,0}{4}{scale=.8} 
\pic{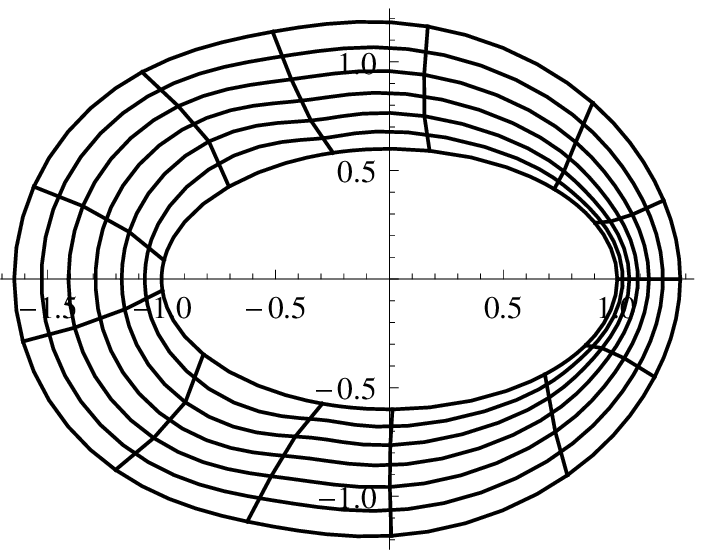}{11,0}{4}{scale=.8}
\caption{Self-mappings of unit disk with Beltrami derivatives
  $\mu_1$ (left), $\mu_2$ (right) followed by conformal mapping to
  exterior of ellipse as in \cite{Dar-fastB}.}
\label{figdar1}
\end{figure}

\example \emph{Exterior mappings.}
In Daripa \cite{Dar-fastB},  quasiconformal mappings
from $\D$ to the exterior of an ellipse (the origin being
sent to $\infty$) are calculated with the
following two sample Beltrami derivatives,
\begin{eqnarray*}
  \mu_1(z) &=& |z|^2 e^{0.65(iz^5-2.0)}, \\
  \mu_2(z) &=& \frac{1}{2}|z|^2 \sin(5\re z) .
\end{eqnarray*}
These exterior mapping results can be related to those of our
algorithm by use of the rational function
\[ h(z) = \frac{(1+\alpha)-(1-\alpha)z^2}{2\alpha z}
\] 
which transforms $\D$ conformally to the exterior
of an ellipse with aspect ratio $\alpha$. 
Composition of $h$ following the quasiconformal self-mapping
of $D$ provides a mapping to the exterior of the ellipse with
the same Beltrami derivative.
 
In the examples in \cite{Dar-fastB}, $\alpha=0.6$ is
specified (however, the inner ellipses in \cite{Dar-fastB} appear
to have aspect ratios of approximately 0.47; axes are not drawn.)
We have made adjustment for the fact that 
Daripa uses $M$ radii equally spaced in $[0,1]$, in contrast to the
exponential spacing we have been using.  
Our results are depicted in Figure \ref{figdar1} with
$(M,N)=(64,64)$.  These images appear fairly similar
those shown in \cite{Dar-fastB}.
   
Computation times are reported in \cite{Dar-fastB} for $N=64$ as
approximately 8.5 seconds of CPU on a MIPS computer described as
``approximately 15 times slower than the CRAY-YMP at Texas A \& M
University'' of that time. Our laptop CPU times, using \emph{Mathematica},  were
approximately 15 seconds for the first example and 45 seconds for the
other. 

\begin{figure}[h!]
\pic{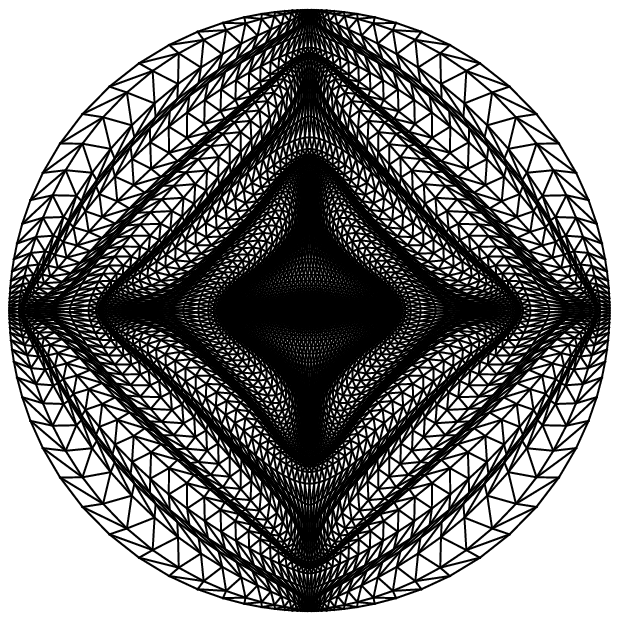}{5,0}{5.5}{scale=.9} 
\pic{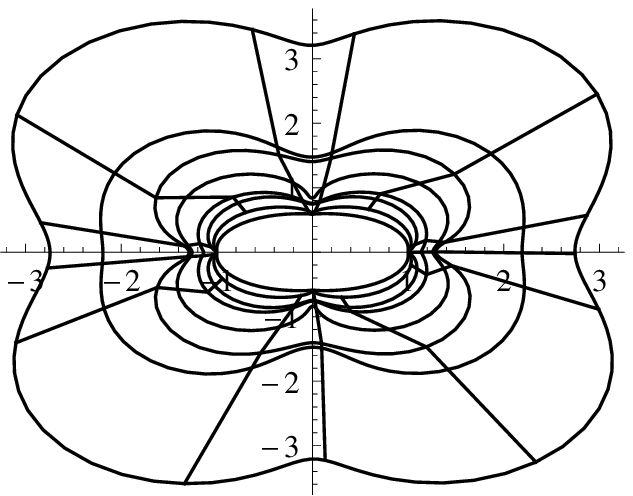}{11,.5}{5}{scale=.9}
\caption{Quasiconformal mappings to disk (left) and
exterior of ellipse (right) determined by (\ref{eq:muoscillate}).}
\label{fig:dartype}
\end{figure}

It is apparent from the examples in \cite{Dar-fastB} that Beltrami
derivatives with a great deal of oscillation were used in order
to create an interesting problem.  Figure \ref{fig:dartype} shows our
results for the mapping of the disk to the exterior of the same ellipse,
with Beltrami derivative
\begin{equation}  \label{eq:muoscillate}
   \mu(z) = 0.9 \sin|20z|
\end{equation}
and $(M,N)=(128,128)$.

\begin{figure}[b!]
  \centering
\pic{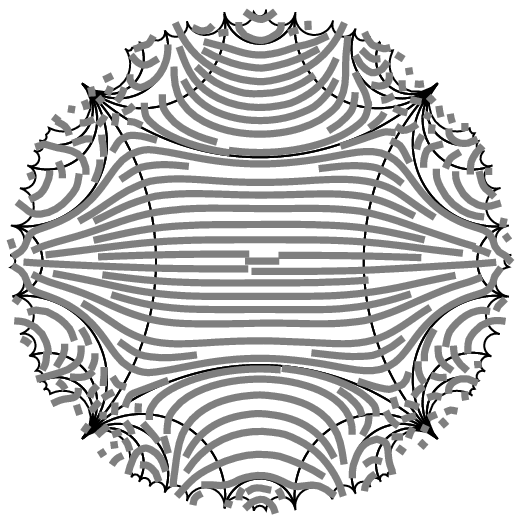}{-.5,-.5}{4.5}{scale=1}
\pic{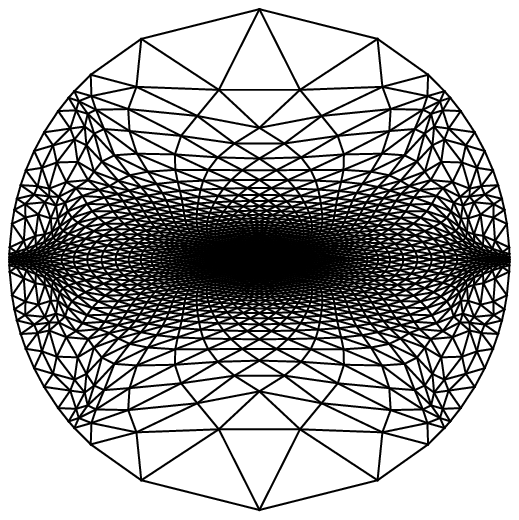}{5.,-.5}{4.5}{scale=1}
\caption{(left) Horizontal trajectories of quadratic differential
  $\Phi$ (\ref{eq:poincareseries}); (right) image triangulation
  arising from induced Teichm\"uller differential
  (\ref{eq:mupoincare}), which effectively produces a streching along
  these trajectories.}\label{fig:torusw}
\end{figure}

\noindent(5) \emph{Quasiconformal deformation of Fuchsian
groups.}  In this example $\mu$ is associated to a quadratic
differential for a Fuchsian group  (see \cite{Harv,Lehto} for
definitions). The two linear-fractional transformations
\[  z\mapsto\frac{z+\sqrt{2}/2}{(\sqrt{2}/2)z +1},\quad
    z\mapsto\frac{z+i\sqrt{2}/2}{-i(\sqrt{2}/2)z +1}
\]
generate a free group $\Gamma$ of self-mappings acting freely on
$\D$. This group has a standard fundamental domain with four-fold
symmetry about the origin as shown in Figure \ref{fig:torusw}, and the
quotient Riemann surface $\D/\Gamma$ is homeomorphic to a torus with a
single puncture.  We create a holomorphic function in $\D$ via the
Poincar\'e series
\begin{equation}\label{eq:poincareseries}
  \Theta(z)=\sum_{\gamma\in G}(\gamma'(z))^2.
\end{equation}
which is known to converge and to satisfy the invariance relation
\begin{equation}\label{eq:poincareinvariances}
  \Theta(\gamma(z))\gamma'(z)^2=\Theta(z)
\end{equation}
for every $\gamma$ in $\Gamma$. Then for $0\le|c|<1$ the function
 \begin{equation} \label{eq:mupoincare}
\mu(z) = c \frac{\conj{\Theta(z)}}{\Theta(z)}
\end{equation}
is a Beltrami differential for the group $\Gamma$; that is,
$\mu(\gamma(z))\conj{\gamma(z)}/\gamma(z) = \mu(z)$
for every $\gamma$ in $\Gamma$.  This implies that the
normalized $\mu$-conformal
self-mapping $f_c\colon\D\to\D$ satisfies $f_c(\gamma(z))=\gamma_c(f_c(z))$
where $\gamma\mapsto\gamma_c$ is an isomorphism from $\Gamma$ to
another Fuchsian group $\Gamma_c$ also acting on $\D$.  A large part
of Teichm\"uller theory is devoted to understanding the nature of
deformed groups such as $\Gamma_c$ and their dependence on $\mu$.
 
\begin{figure}[t!]
\pic{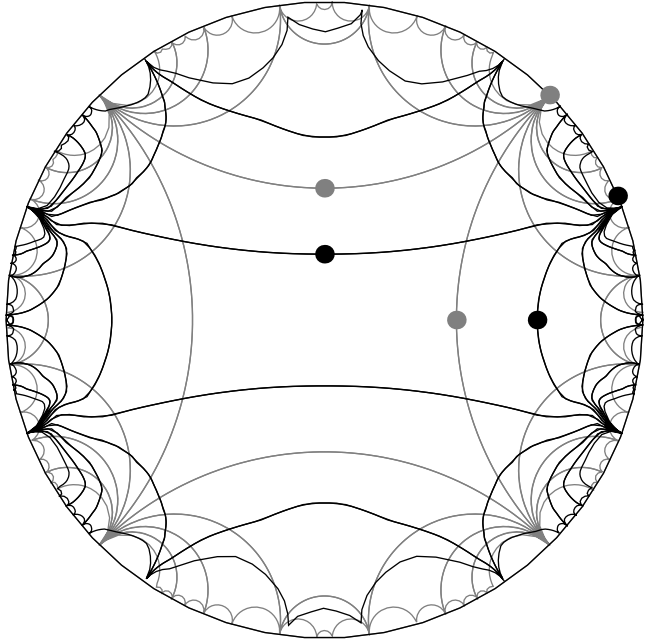}{6.5,0}{8}{scale=1}
\caption{Self-mapping of unit disk, with Beltrami derivative defined
  via a Poincar\'e series.  Computations were made with
  $(M,N)=(64,64)$ The standard fundamental domain (lighter contour) is
  symmetric under rotation by $\pi/2$; its image under the
  quasiconformal mapping (black contour) is superimposed.}
\label{fig:torus1}
\end{figure}

Figure \ref{fig:torusw} shows the trajectory structure of the
quadratic differential $\Phi$, that is, the solutions of
$\Phi(z)dz^2>0$, together with a $w$-triangulation obtained from the
$\mu$-conformal mapping for $c=0.5$. Figure \ref{fig:torusw} shows the
$w$-triangulation for $(M,N)=(64,64)$.  Figure \ref{fig:torus1} is a
superposition of the $z$- and $w$-domains, with part of the
tesselation of $\D$ by group translates of the fundamental domain of
$\Gamma$ and also the tessellation by their images, which are
fundamental domains for $\Gamma_\mu$.  The images were calculated by
taking each point defining each curve of the original tessellation,
identifying the particular $z$-triangle in which it lies, and then
applying the PL-mapping to the corresponding $w$-triangle of Figure
\ref{fig:torusw}.  Some inaccuracies, particularly near the boundary,
are clearly visible inasmuch as the image curves must be hyperbolic
geodesics. One probable source of error is in the calculation of
$\Theta(z)$. To evaluate (\ref{eq:poincareseries}) we truncated the
series to words of up to length 6 in the generators of $\Gamma$ and/or
inverses, and used (\ref{eq:poincareinvariances}) to apply
(\ref{eq:poincareseries}) directly only for $z$ inside the fundamental
domain containing the origin; then (\ref{eq:mupoincare}) was applied
for points outside the fundamental domain.

\begin{figure}[t!]
\pic{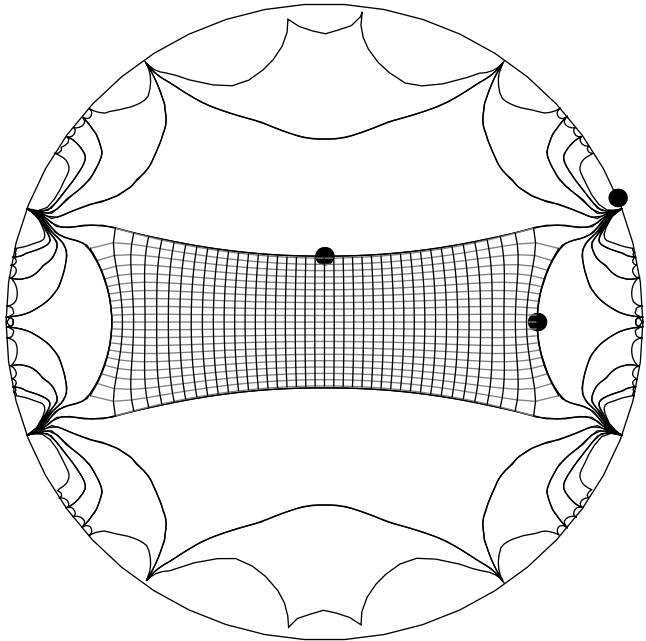}{6.5,0}{8}{scale=1}
\caption{Comparison of the image of the basic fundamental
domain via quasiconformal mapping and conformal mapping techniques.
The  trajectories inside the fundamental domain were
calculated as the images under the conformal mapping of vertical
and horizontal segments in the period rectangle of the Weierstrass
$\wp$-function.}
\label{fig:torus2}
\end{figure}

An independent verification of this mapping is made using the results of
\cite{BP} on conformal mapping of symmetric quadrilaterals with
circular sides.  In the present context, where the vertex angles of
the fundamental domains are all equal to zero, such quadrilaterals
depend (up to conformal equivalence) on two real parameters, which may
be prescribed in several ways, some relating to the geometry of the
circular quadrilateral (such as the midpoints $p_1$, $p_2$ of the right and
upper edges as marked in Figure \ref{fig:torus1}), others relating to
parameters in the conformal mapping.  In particular, there is a
conformal mappping $h$ from a Euclidean rectangle with vertices at
$\pm\omega_1\pm\omega_2$ ($\omega_1>0$, $\omega_2/i>0$) to the
fundamental domain quadrilateral, whose Schwarzian derivative $\mathcal{S}_h$
\cite{AhlLQM} is equal to the following elliptic function:
\begin{equation}  \label{eq:schwarzian}
  \mathcal{S}_h(\zeta)=\frac{1}{2}\wp(\zeta+\omega_1+\omega_2) - 2\sigma,
\end{equation}
where the essential parameter is the purely imaginary ratio
$\tau=\omega_2/\omega_1$, which determines the second parameter
$\sigma=\sigma(\tau)\in\R$ due to the fact that the circular edges
must be orthogonal to the circle passing through the four image
vertices $h(\pm\omega_1\pm\omega_2)$.  Here $\omega_1,\omega_2$ are
basic half-periods of the Weierstrass $\wp$-function, normalized in a
specific way which we will not describe here. Methods are given in
\cite{BP} for calculating the relationships among these parameters. We
have applied them to the observed values for $p_1=f_c(\sqrt{2}-1)$ and
$p_2=f_c(i(\sqrt{2}-1))$ obtained from our algorithm for the Beltrami
equation, together with the fact that $\tau=(1-c)/(1+c)$, to obtain
the value of $\sigma$ numerically.  Then the Schwarzian differential
equation (\ref{eq:schwarzian}) was solved, first along the real and
imaginary axes, and from there along vertical and horizontal lines
throughout the rectangle with vertices $\pm\omega_1\pm\omega_2$.
Since this was done with a normalization of $h'(0)=1$, a further
scaling was necessary to adjust the size to match the value of $p_1$.
It is seen (Figure \ref{fig:torus2}) that the midpoint $p_2$ and the
curvatures of the edges agree quite well.  In general if one is
interested in deforming Fuchsian or Kleinian groups by this method, a
reasonable strategy would be to replicate the information obtained in
such a fundamental domain via the group action of $\Gamma_c$.
  
\begin{table}[!th]  \centering
$ \begin{array}{l|l||cccc}
 \multicolumn{1}{c|}{c} & 
 \multicolumn{1}{c||}{\tau} & 
 N=16 & N=32 & N=48 & N=64  \\\hline
 0.1 & 0.818 & 0.022 & 0.006 & 0.003 & 0.002\\
 0.3 & 0.538 & 0.089 & 0.027 & 0.011 & 0.008\\
 0.5 & 0.333 & 0.283 & 0.110 & 0.052 & 0.031\\
 \end{array}$
 \caption{Relative error of midpoint $p_2$ of the upper edge
of the deformed fundamental domain compared to result of conformal mapping from rectangle. Calculations made with $M=N$.}
  \label{tab:torus}
\end{table}

\section{Discussion and Conclusions}

It is stated in Daripa \cite{Dar-fastB} that prior to that article
there were no constructive methods published for solving the Beltrami
equation in the disk numerically. Convergence proofs appeared later in
\cite{GK} together with a modified scheme. There are other methods
which also include convergence proofs (we have mentioned \cite{He,W})
without giving a detailed analysis of the rate of convergence.
Therefore we discuss here some aspects of \cite{Dar-fastB,GK} in
relation to our algorithm.  As we mentioned in the Introduction, that
approach is based on evaluation of singular integrals. The original
problem is presented in the context of finding a $\mu$-conformal
mapping to a prescribed star-shaped domain, and is in some ways
reminiscent of the classical method of Theodorsen \cite{Henr} for
conformal mappings.

Daripa's main algorithm requires evaluation of the
$\partial/\partial\overline{z}$ derivatives which appear in the
singular integrals.  A variant is also proposed which does not require
these derivatives; however this is not applied in the numerical
examples provided.  The operation count of one iteration of Daripa's
method is $O(MN\log N)$. This should be multiplied by the average
number of iterations required, which depends on how refined the mesh
is and how much accuracy is desired. In the examples which we have
taken from \cite{Dar-fastB}, $\|\mu\|_\infty$ is approximately $0.5$,
but it should be noted that $|\mu(z)|$ is in fact bounded by 0.12 for
$|z|<0.5$, and by 0.05 for $|z|<0.3$.  In fact, an important
limitation stated in \cite{Dar-fastB} is that the Beltrami derivative
$\mu$ must be H\"older continuous. Further, it is recommended that
$\mu$ vanish at least as fast as $|z|^3$ at the origin for the method
to work properly. In \cite{GK} it is similarly recognized that
computation time increases as $\|\mu\|_\infty$ increases.  Our
algorithm, in contrast, is not subject to any such requirement on
$\mu$. Our computation times are considerably longer than those
reported in \cite{GK} but this may be due in large part to use of a
symbolic interpreter rather than a compiled program.

We also mention the method presented in \cite{LLSX} for determining a
Teichm\"uller mapping to an arbitrary domain. It begins by choosing an
``optimal'' Beltrami derivative and then solving the corresponding
boundary Beltrami equation.  A system of linear equations is
determined to discretize this equation, which must be solved together
with a collection of nonlinear boundary constraints.  This is solved
by an iterative method (conjugate gradient).

Our algorithm involves no evaluation of singular integrals and no
iteration of solutions.  As described in Section \ref{sec:main}, we
use a purely linear system $(A,B)$.  For a mesh of dimensions $M,N$,
the matrices $A,B$ are of orders $n_e \times n_{\rm v}$, $n_{\rm e}\times1$
respectively, where $n_{\rm e}=4MN+2(N-1)+1$ and $n_{\rm v}=(2M+1)N$.  Although
the total number of elements contained in $A$ is $O(M^2N^2)$, by
construction $A$ is a sparse matrix: as we noted during the proof of Theorem
\ref{th:main}, the number of nonzero elements of
any row of $A$ is no greater than 3; further, the
number of nonzero elements of any column of $A$ is no greater than 7
(note that each variable $W_{jk}$ in (\ref{eq:triangles+}), (\ref{eq:triangles-})
corresponds to a vertex of at most six triangles, cf. Figure
\ref{fig:Wmesh}; the only vertex appearing in seven equations is $W_{0,0}$, cf.\
(\ref{eq:normalize})). Hence the number of nonzero elements in $A$ is no
greater than $O(MN)$. 

Since the linear system $AV=B$ is highly overdetermined and in general
there is no exact solution, we use the least squares approximation.
There are many numerical methods available for the least squares
problem; a comprehensive reference is \cite{Bj}.  For simplicity we
discuss the method of ``normal equations,'' that is, the solution of
$A^H A V =A^H B$, where $A^H$ is the conjugate transpose of $A$.  It
is easily seen that by construction the columns of $A$ are linearly
independent, so that $A^H A$ is positive definite, and each row or
column of $A^H A$ contains at most $7$ nonzero entries.  
These entries are not consecutive, but the index correspondence in
(\ref{eq:VpWjk}) may be taken so that the row bandwidth of $A^H A$ is
$2M$. When one applies Gaussian reduction to $(A,B)$, the
$n$th of the $n_{\rm v}$ rows will only need to reduce at most $2M$ of
the succeeding rows, and will require no more than $2M$ floating-point
multiplications for each one.  Thus the total operation count is of
the order of $O(n_{\rm v}(2M)^2)=O(M^3N)$.  Once $A^H A$ has been thus
reduced to echelon form, the computational cost of back substitution
is seen to be no more than $O(M^2 N)$. These computations do not
require a significant amount of storage other than the original
data. In summary, when one doubles the mesh dimensions $M,N$, the
memory requirement is at most multiplied by $4$ (which is the same as
the increase in the mesh itself) and the computation time by 16.  Our
numerical experiments indicate that the more sophisticated
least-squares algorithms found in packaged software appear to reduce
these exponents slightly.
 
It may be noted that our algorithm can be used for solving the
Beltrami equation on the entire $z$-plane (normalized by fixing 0, 1,
and $\infty$) instead of the disk.  One simply eliminates the step of extending the
Beltrami coefficient from the disk to its exterior by reflection.  We
have not yet investigated this question numerically.  Further, it
seems probable that our algorithm will converge to the $\mu$-conformal
mapping even when $\mu$ is only piecewise smooth, and perhaps in even
greater generality. Many numerical examples suggest this, some of
which we have given above. The estimates in the proof
of Theorem \ref{th:main} concerning $L_2$ norms would not be greatly
affected if only a small proportion of the terms in the sums failed to
tend to zero as fast as required.  We hope to look into these
questions in future work.

 We believe that this algorithm for solving the Beltrami equation is
conceptually much simpler than other methods which have been presented, and
is easy to implement.

The authors are grateful to T.\ Sugawa for many critical and useful
comments in the preparation of this work. The first author is also
grateful to D.\ Marshall for pointing out several fundamental errors
in the approach proposed in \cite{PorArx}.


\begin{thebibliography}{99}
 
\bibitem{AhlLQM} L.\ Ahlfors, \emph{Lectures on Quasiconformal
Mappings}, second edition, University Lecture Series \textbf{38},
American Mathematical Society, Providence, RI (2006).

\bibitem{AHTK} S. Angenent, S. Haker, A. Tannenbaum, and
R. Kikinis, Laplace-Beltrami operator and brain surface flattening,
\emph{IEEE Trans.\ Medical Imaging} \textbf{18}:700--711 (1999).


\bibitem{AMPP} K.\ Astala, J.\ L.\ Mueller, A.\ Per\"am\"aki, L.\
  P\"aiv\"arinta, and S. Siltanen. Direct electrical impedance
  tomography for nonsmooth conductivities, \emph{Inverse Problems and
    Imaging} \textbf{5}:531--549 (2011).

\bibitem{Bj} \AA. Bj\"orck, \emph{Numerical Methods for Least
    Squares Problems,} SIAM. (1996) ISBN 978-0-89871-360-2.

 
\bibitem{BP} P. Brown and R. M. Porter, Conformal mapping of circular
  quadrilaterals and Weierstrass elliptic functions, \emph{Comput.\
    Methods Funct.\ Theory} \textbf{11}:463--486 (2011).

\bibitem{Dar-Mash} P.~Daripa and D.\ Mashat, An efficient and novel
numerical method for quasiconformal domains,
\emph{Numer. Algorithms} \textbf{18}:159-175 (1998).

\bibitem{Dar-fastB} P.\ Daripa, A fast algorithm to solve
the Beltrami equation with applications to quasiconformal
mappings,  \emph{J.\ Comput.\ Phys.} \textbf{106}:355-365 (1993).

\bibitem{Dar-ap} P.\ Daripa, On applications of a complex variable 
method in compressible flows, \emph{J.\ Comput.\ Phys.}
\textbf{88}:337--361 (1990).
 
\bibitem{Dar-CR} P.\ Daripa, A fast algorithm to solve
nonhomogeneous Cauchy-Riemann equations in the complex plane,
\emph{SIAM J.\ Sci.\ Statist.\ Comput.} \textbf{13}:1418--1432
(1992).


 
\bibitem{DT} T.\ A.\ Driscoll and L.\ N.\ Trefethen
\emph{Schwarz-Christoffel mapping}, Cambridge Monographs on
Applied and Computational Mathematics \textbf{8}, Cambridge University
Press, Cambridge (2002). 

\bibitem{GK} D.\ Gaidashev, D.\ Khmelev, On numerical algorithms for the solution of a Beltrami equation,. \emph{SIAM J.\ Numer.\ Anal.}\ \textbf{46}:5 (2008) 2238--2253.  



 

\bibitem{GZWLY} X.\ D.\ Gu, W.\ Zeng, F.\ Luo, and Sh-T.\ Yau,
  Numerical computation of surface conformal mappings,
  \emph{Comput.\ Methods Funct.\ Theory} \textbf{11}:747--787(2011).

 
\bibitem{Harv} W.\ J.\, Harvey, ed.\ \emph{Discrete Groups and
Automorphic Functions.  Proceedings of an Instructional Conference
held in Cambridge, July 28--August 15, 1975}, Academic Press 
(1977).

\bibitem{He} Zh.-X.\ He, Solving Beltrami equations by circle packing,
\emph{Trans.\ A.M.S.}\  \textbf{322}:657--670 (1990). 

\bibitem{Henr} P.\ Henrici, \emph{Applied and Computational Complex
Analysis,} vol.\ \textbf{3}, Wiley, New York (1986).


\bibitem{KaTs2} D.\ S.\ Kamenetski\u\i\ and S.\ V.\ Tsynkov, On the
construction of images of simply connected domains realized by
solutions of a system of Beltrami equations (Russian),
\emph{Akad.\ Nauk SSSR Inst.\ Prikl.\ Mat.}  Preprint no. 155 (1990).
    

\bibitem{Lehto} O.\ Lehto, \emph{Univalent Functions and
Teichm\"uller Spaces}, Springer-Verlag, New York (1995).

\bibitem{LV} O.\ Lehto and K.\ I.\ Virtanen, \emph{Quasiconformal
Mappings in the Plane}, second edition, \emph{Die Grundlehren der
mathematischen Wissenschaften} \textbf{126} Springer-Verlag, New
York-Heidelberg (1973).

\bibitem{LPRM} B.\ L\'evy, S.\ Petitjean, N.\ Ray and J.\ Maillot,
  Least Squares Conformal Maps for Automatic Texture Atlas
  Generation, \emph{ACM Transactions on Graphics (TOG)}, Proceedings
  of ACM SIGGRAPH 2002 \textbf{21}:362--371 (2002).

\bibitem{LLKWG} L.\ M.\ Lui, K.\ C.\ Lam, T.\ W. Wong, and X.\ Gu,
Texture Map and Video Compression Using Beltrami Representation,
\emph{SIAM J.\ Imaging Sci.} \textbf{6}:1880--1902 (2013).

\bibitem{LLSX} L.\ M.\ Lui, K.\ C.\ Lam, S-T.\ Yau, and X.\ Gu,
Teichm\"uller extre\-mal mapping and its applications to landmark
matching registration, \texttt{arXiv:1211.2569v1}.


\bibitem{Mu} J. R. Munkres, \emph{Elementary Differential Topology. Lectures
  given at Massachusetts Institute of Technology, Fall, 1961}. Revised
  edition. Annals of Mathematics Studies, No. 54 Princeton University
  Press, Princeton, N.J.




\bibitem{PorArx} R.\ M.\ Porter,  Numerical solution of the Beltrami equation,
\texttt{arXiv: 0802.1195} [math.CV]


\bibitem{Sam} Z.\ V.\ Samsoniya, Construction of certain
quasiconformal mappings (Russian), \emph{Trudy Vychisl.\ Tsentra
Akad.\ Nauk Gruzin.\ SSR} \textbf{23}:76-89 (1983).    


\bibitem{Sz} G.\ Szeg\"o, Conformal mapping of the interior of
an ellipse onto a circle, \emph{American Mathematical Monthly}
\textbf{57}:474--479 (1950).

\bibitem{Tr} M.\ Trott, \emph{The Mathematica GuideBook for Numerics},
Springer Science+Business Media, Inc. (2006).



\bibitem{WW} E.\ T.\ Whittaker and G.\ N.\ Watson, \emph{A Course
of Modern Analysis}, reprint of the fourth (1927) edition,
Cambridge Mathematical Library, Cambridge University Press,
Cambridge (1996). 

\bibitem{W} G.\ B.\ Williams, A circle packing measurable
Riemann theorem, \emph{Proc.\ Amer.\ Math.\ Soc.} 
\textbf{134}:2139--2146 (2006). 



\end{thebibliography}
\end{document}